\documentclass[a4paper,12pt]{article}
\usepackage{amsfonts}
\usepackage{color}
\usepackage{eurosym}
\usepackage[utf8]{inputenc}
\usepackage[english]{babel}
\usepackage{amsmath,amsxtra,amssymb,latexsym, amscd,amsthm}
\usepackage[mathscr]{eucal}
\usepackage{latexsym, amsmath, amsfonts, amscd, amssymb, verbatim,color,cite}
\def\tto{\;{\lower 1pt \hbox{$\rightarrow$}}\kern -10pt
\hbox{\raise 2pt \hbox{$\rightarrow$}}\;}

\def\B{I\!\!B}

\def\R{I\!\!R}
\def\N{I\!\!N}

\def \N{I\!\!N}

\begin{document}
\pagestyle{myheadings}
\newtheorem{Theorem}{Theorem}[section]
\newtheorem{Proposition}[Theorem]{Proposition}
\newtheorem{Remark}[Theorem]{Remark}
\newtheorem{Lemma}[Theorem]{Lemma}
\newtheorem{Corollary}[Theorem]{Corollary}
\newtheorem{Definition}[Theorem]{Definition}
\newtheorem{Example}[Theorem]{Example}
\renewcommand{\theequation}{\thesection.\arabic{equation}}
\normalsize

\title{\sf Robust Global  Solutions of  Bilevel  Polynomial Optimization Problems with   Uncertain Linear Constraints}

\author{\sc T. D. Chuong\footnote{School of Mathematics and Statistics,
University of New South Wales, Sydney NSW 2052, Australia; email:
chuongthaidoan@yahoo.com. Research was
supported by the UNSW VC's Postdoctoral fellowship program.} \  and \
 V. Jeyakumar\footnote{School of Mathematics and Statistics, University of New South Wales,
Sydney NSW 2052, Australia; email:
 v.jeyakumar@unsw.edu.au. Research was partially supported by a grant from the
Australian Research Council}}


 \maketitle

\begin{abstract}
This paper studies, for the first time, a bilevel polynomial program whose constraints involve uncertain linear constraints and another uncertain linear optimization problem. In the case of box data uncertainty, we present a sum of squares polynomial characterization of a global solution of its robust counterpart where the constraints are enforced  for all realizations of the uncertainties within the prescribed uncertainty sets. By characterizing a solution of the robust counterpart of the lower-level uncertain linear program under spectrahedral uncertainty using a new generalization of Farkas' lemma, we reformulate the robust bilevel program as a single level non-convex polynomial optimization problem. We then characterize a global solution of the single level polynomial program by employing Putinar's Positivstellensatz of algebraic geometry under coercivity of the polynomial objective function. Consequently, we show that the robust global optimal value of the bilevel program is the limit of a sequence of values of Lasserre-type hierarchy of semidefinite linear programming relaxations. Numerical examples are given to show how the robust optimal value of the bilevel program can be calculated by solving semidefinite programming problems using the Matlab toolbox YALMIP.

\medskip
{\bf Key words.} Bilevel Programming, uncertain linear constraints, robust optimization, global polynomial optimization.

\end{abstract}
\markboth{\centerline{\sc }}{\centerline{\sc {t. d. chuong and v. jeyakumar}}}
\section{Introduction}
\setcounter{equation}{0}

The bilevel optimization problems arise when two independent decision makers, ordered within a hierarchical structure, have conflicting objectives. They appear as hierarchical decision-making problems, such as risk management and economic planning problems, in engineering, governments and industries \cite{Bilevelbook,dempeBook,Marcotte,Dempe0,Floudas2001}. The commonly used  bilevel optimization techniques (see \cite{Dempe0,Dempe_Z,ye1,ye2} and other references therein) assume perfect information (that is, accurate values for the input
quantities or system parameters), despite the reality that such precise knowledge is rarely available in hierarchical decision-making problems. The data of these problems are often uncertain (that is, they are not known exactly at the time of the decision) due to estimation errors, prediction errors or lack of information. Consequently, the development of  optimization methodologies which are capable of generating robust optimal solutions that are immunized against data uncertainty,  such as the deterministic robust optimization techniques, has become more important than ever in mathematics, commerce and engineering \cite{robustbook,siamRevRobust,Goberna-Jeya-Li15}. Yet, such a mathematical theory and the associated methods for bilevel optimization in the face of data uncertainty do not appear to be available in the literature.
\medskip

In this paper  we study, for the first time, the following bilevel polynomial program  that finds robust global optimal solutions under bounded data uncertainty:
{\small
\begin{eqnarray*}
 &\displaystyle\min_{(x,y)\in \R^{m}\times \R^{n}}&
f(x,y) \\
&\mbox{ subject to }&  \tilde a_i^\top x+\tilde b_i^\top y\le \tilde c_i,\; \forall (\tilde a_i, \tilde b_i, \tilde c_i)\in\widetilde U_i,\; i=1,\ldots,l, \\
&&y\in Y(x):={\rm argmin}_{z\in\R^n}\{c_0^\top x+d_0^\top z\mid c_j^\top x +\hat a_j^\top z\le \hat b_j,\; \forall (\hat a_j,\hat b_j)\in \widehat{U}_j,\;  j=1,\ldots,q\}.
\end{eqnarray*}}
It is the robust counterpart of the bilevel polynomial program with uncertain linear constraints
{\small
\begin{eqnarray*}
 &\displaystyle\min_{(x,y)\in \R^{m}\times \R^{n}}&
f(x,y) \\
&\mbox{ subject to }&  \tilde a_i^\top x+\tilde b_i^\top y\le \tilde c_i, \; i=1,\ldots,l, \\
&&y\in {\rm argmin}_{z\in\R^n}\{c_0^\top x+d_0^\top z\mid c_j^\top x +\hat a_j^\top z\le \hat b_j,\;  j=1,\ldots,q\},
\end{eqnarray*}}
\vspace{-0.8cm}

\noindent where the constraint data $(\tilde a_i, \tilde b_i, \tilde c_i), i=1,\ldots,l$ and $(\hat a_j,\hat b_j), j=1,\ldots, q$ are {\it uncertain} and they belong to the prescribed bounded {\it uncertainty} sets $ \widetilde U_i, i=1,\ldots,l$ and $\widehat U_j, j=1,\ldots, q$, respectively. The vectors $c_0\in\R^m, d_0\in\R^n, c_j \in \R^m,    j=1,\ldots,q$ are fixed and $f(x,y)$ is a polynomial. Note that, in the robust counterpart, the uncertain linear constraints are enforced for all realizations of the uncertainties within the uncertainty sets. The robust counterpart is, therefore, a worst-case formulation in terms of deviations of the data from their nominal values.
\medskip


The bilevel program is a class of hard optimization problems even for the case where all the functions are linear and are free of data uncertainty \cite{hard}.
A general approach for studying bilevel optimization problems is to transform them into  single level optimization problems \cite{Marcotte,Dempe0,Dempe_Z}. The resulting single level optimization problems are generally non-convex constrained optimization problems. It is often difficult to find global optimal solutions of non-convex optimization problems. However, Putinar Positivstellensatz \cite{Putinar93} together with Lasserre type semidefinite relaxations allows us to characterize global optimal solutions and find the global optimal value of a non-convex optimization involving polynomials. The reader is referred to \cite{JLL,Lasserre-book-09,jeya-Son-Li-ORL14,jeya-optimL,JKLL,JLLP} for related recent work on single level convex and non-convex polynomial optimization in the literature.
\medskip

In this paper we make the following key contributions to bilevel optimization.
\medskip

\vspace{-0.6cm}
\begin{itemize}
 \item{\bf Characterizations of lower-level robust solutions \& spectrahedral uncertainty.}
We first present a complete dual characterization for a robust solution of the lower-level uncertain linear program in the general case of bounded spectrahedral uncertainty sets $U_j, j=1,2,\ldots,q$. It is achieved by way of proving a generalization of the celebrated non-homogeneous Farkas' lemma \cite{dinh-jeya,Hiriart-Urruty93} for semi-infinite linear inequality systems whose dual statement can be verified by solving a semidefinite linear program. A special variable transformation paves the way to formulate the dual statement in terms of linear matrix inequalities.  The spectrahedral uncertainty set possesses a broad spectrum of convex uncertainty sets that appear in robust optimization. It includes polyhedra,  balls and ellipsoids\cite{Ramana-Goldman95,Vinzant14}, for which dual characterizations of robust optimality of uncertain linear programs are already available in the literature  \cite{Ben-tal01,robustbook,siamRevRobust}.


 \item{\bf A sum of squares polynomial characterization of robust global optimal solutions.}
 In the case of box data uncertainty in both upper and lower level constraints and the objective polynomial is coercive, we  derive a sum of squares polynomial characterization of robust global solution of the bilevel program by first transforming the bilevel program into a single level non-convex polynomial program using the dual characterization of solution of the lower level program and then employing the Putinar Positivstellensatz \cite{Putinar93}. It is derived under a suitable Slater-type regularity condition. A numerical example is given to show that our characterization may fail without the regularity condition. In the Appendix, we show how a
numerically checkable characterization of robust feasible solution can be obtained in the case of ball data uncertainty in the constraints. Related recent work on global bilevel polynomial optimization in the absence of data uncertainty can be found in \cite{Jeya-Li-bilevel15,JLLP}.


\item{\bf Convergence of SDP relaxations to the robust global optimal value.}
Finally, using our sum of squares characterization of robust global optimality together with Lasserre type semidefinite relaxations, we show how robust global optimal value can be calculated by solving a sequence of semidefinite linear programming relaxations. We prove that the values of the Lasserre type semidefinite relaxations converge to the global optimal value of the bilevel polynomial problem. We provide numerical examples to illustrate how the optimal value can be found using the Matlab toolbox YALMIP.
\end{itemize}

The outline of the paper is as follows. Section 1 presents a generalization of the Farkas lemma and a dual characterization for robust global optimality of the lower level program of  the bilevel  problem. Section 3 develops characterizations for robust solution of uncertain bilevel problems  in the case of box data uncertainty in the constraints. Section 4 provides results on finding the optimal value of  the bilevel  problem via semidefinite programming relaxations. Appendix presents numerically checkable characterizations of robust feasible solutions of the bilevel program  in the case of ball data uncertainty.

\section{Lower Level LPs under Spectrahedral Uncertainty}
\setcounter{equation}{0}
In this section, we present a robust non-homogeneous Farkas's lemma for an uncertain linear inequality system.

\medskip
Let us first recall   some  notation and preliminaries. The notation $\R^n$ signifies  the Euclidean  space whose norm is denoted by
$\|\cdot\|$ for each $n\in\N:=\{1,2,\ldots\}$. An element $x\in \R^n$ is written as a column vector, but it is sometimes convenient to write the components of a vector in a row instead of a column. The inner product in $\R^n$ is defined by $\langle x,y\rangle:=x^\top y$ for all $x, y\in\R^n.$  The topological   closure  of a set $\Omega\subset \R^n$ is
denoted by ${\rm cl\,}\Omega$.   The origin of any space is denoted by $0$ but  we may use
$0_n$  for the origin of $\R^n$ in situations where some confusion might be
possible. As usual, ${\rm conv\,}\Omega$ denotes the convex hull of $\Omega,$ while ${\rm
cone\,}\Omega:=\R_+{\rm conv\,}\Omega$ stands for the convex conical hull of
$\Omega\cup \{0_n\}$, where $\R_+:=[0,+\infty)\subset \R.$ A symmetric $(n\times n)$  matrix $A$ is said to be {\it positive semi-definite}, denoted by $A\succeq 0$, whenever $x^\top Ax\ge 0$ for all $x\in\R^n.$


\medskip
The classical version of Farkas's Lemma for a semi-infinite inequality system can be found in \cite[Theorem~4.3.4]{Hiriart-Urruty93}.

\begin{Lemma}\label{Ge-Farkas} {\rm (Non-homogeneous Farkas' lemma--\cite[Theorem~4.3.4]{Hiriart-Urruty93})} Let $(b,r), (s_j,p_j)\in \R^n\times\R$, where $j$ varies in an arbitrary index set $J$. Suppose that the system of inequalities \begin{align}\label{G-F} s_j^\top x\le p_j\mbox{ for all } j\in J\end{align} has a solution $x\in\R^n.$ Then, the following two properties are equivalent:

{\rm (i) }  $b^\top x\le r \mbox{ for all $x$ satisfying } \eqref{G-F};$

{\rm (ii) } $(b,r)\in {\rm cl\,cone}\{(0_n,1)\cup(s_j,p_j)\mid j\in J\}.$
\end{Lemma}

Consider an {\it uncertain} linear inequality system,
\begin{align}\label{un-orginal}x\in\R^n,\; \hat a_j^\top x\le \hat b_j,\; j=1,\ldots, q,\end{align}
where  $(\hat a_j, \hat b_j), j=1,\ldots, q$ are uncertain and they belong to the {\it uncertainty} sets $\widehat U_j, j=1,\ldots, q.$ The uncertainty sets are given by
$$\widehat U_j:=\Big\{(a^0_j+\sum^{s}_{i=1}u^i_ja^i_j,b^0_j+\sum^{s}_{i=1}u^i_jb^i_j)\mid  u_j:=(u^1_j,\ldots, u^s_j)\in U_j\Big\}, \;j=1,\ldots,q,$$ where $a_j^i\in \R^n, b_j^i\in\R, i=0,1,\ldots,s, j=1,\ldots,q$ are fixed and $U_j$ is a {\it spectrahedron} \cite{Ramana-Goldman95,Vinzant14} described by   \begin{align}\label{U-set}U_j:=\big\{u_j:=(u^1_j,\ldots, u^s_j)\in\R^s\mid A^0_j+\sum_{i=1}^su^i_jA^i_j\succeq 0\big\},\; j=1,\ldots, q,\end{align}  and $A^i_j, i=0,1,\ldots,s, j=1,\ldots,q,$ are symmetric $(p\times p)$  matrices.

\medskip
Now, consider the affine mappings  $a_j:\R^s\to \R^n$ and $b_j:\R^s\to\R, j=1,\ldots,q$ given  respectively by \begin{align}\label{affine-maps} a_j(u_j):=a^0_j+\sum^{s}_{i=1}u^i_ja^i_j,\; b_j(u_j):=b^0_j+\sum^{s}_{i=1}u^i_jb^i_j \; \mbox{ for } \; u_j:=(u^1_j,\ldots, u^s_j)\in \R^s\end{align}  with $a_j^i\in \R^n, b_j^i\in\R, i=0,1,\ldots,s, j=1,\ldots,q$ fixed as above. Then, the robust counterpart of the uncertain system~\eqref{un-orginal},
\begin{align*}\label{}x\in\R^n,\; \hat a_j^\top x\le \hat b_j,\;  \forall (\hat a_j,\hat b_j)\in \widehat U_j,\; j=1,\ldots,q,\end{align*}
can be expressed equivalently as \begin{align}\label{un-sys}x\in\R^n,\; a_j(u_j)^\top x\le b_j(u_j),\;  \forall u_j\in U_j,\; j=1,\ldots,q.\end{align}

From now on, the sets $U_j, j=1,\ldots,p$ in \eqref{U-set} are assumed to be  {\it compact}.  Note that the spectrahedra \eqref{U-set} are closed and convex sets, and they possess a broad spectrum of infinite convex sets,  such as polyhedra,  balls, ellipsoids and cylinders \cite{Ramana-Goldman95,Vinzant14}, that appear in robust optimization.

\medskip
In the next theorem, by employing a variable transformation together with Lemma 2.1, we derive a generalized non-homogeneous Farkas's lemma, which provides a numerically tractable certificate for nonnegativity of an affine function over  the uncertain linear inequality system~\eqref{un-sys}. This nonnegativity can be checked by solving a feasibility problem of a semi-definite linear program and it plays an important  role in characterizing robust solutions of the bilevel polynomial optimization problem. For recent work on generalized Farkas' lemma for uncertain linear inequality systems, see \cite{dinh-jeya,goberna-jeya-MP,goberna-jeya-siam}.

\begin{Theorem}\label{Theo1}{\bf (Generalized Farkas' Lemma with numerically checkable dual condition)} Let $(\wp,r)\in\R^n\times\R$.  Assume that the cone $C:={\rm cone}\big\{\big(a_j(u_j),b_j(u_j)\big)\mid u_j\in U_j, j=1,\ldots,q\big\}$ is closed. Then, the following statements are equivalent:

{\rm (i)\ } $x\in\R^n,\; a_j(u_j)^\top x\le b_j(u_j),\;  \forall u_j\in U_j,\; j=1,\ldots,q\Longrightarrow \wp^\top x -r \ge 0;$

{\rm (ii)\ } $\exists  \lambda_j^0\ge 0, \lambda_j^i\in\R, j=1,\ldots,q, i=1,\ldots,s \mbox{ such that}$ \begin{align*} &\wp+\sum\limits_{j=1}^q(\lambda_j^0a^0_j+\sum^{s}_{i=1}\lambda^i_ja^i_j)=0,\; -r-\sum\limits_{j=1}^q(\lambda_j^0b^0_j+\sum^{s}_{i=1}\lambda^i_jb^i_j)\ge 0\\& \mbox{and } \lambda_j^0A^0_j+\sum_{i=1}^s\lambda^i_jA^i_j\succeq 0,  j=1,\ldots,q.\end{align*}
\end{Theorem}
{\bf Proof.} [{\bf (i) $\Longrightarrow$ (ii)}]  Assume that (i) holds. Putting  \begin{align}\label{coneC-tilde}\tilde C:={\rm cone}\Big\{(0_n,1)\cup \big(a_j(u_j),b_j(u_j)\big)\mid u_j\in U_j, j=1,\ldots,q\Big\},\end{align} we first prove that $\tilde C$ is closed. Consider a sequence $x^*_k\to x^*,$ where $x^*_k\in \tilde C,  k\in\N$. Then, $x^*_k:=\lambda_0^k(0_n,1)+\sum\limits_{j=1}^q\lambda_j^k\big(a_j(u_j^k),b_j(u_j^k)\big)$ with $\lambda_j^k\ge 0, j=0,1,\ldots,q$ and $u_j^k\in U_j,  j=1,\ldots,q.$ We assert that the sequence $\lambda_0^k, k\in \N$ is bounded. Otherwise, by taking a subsequence if necessary we may assume that $\lambda_0^k\to+\infty$ as $k\to\infty.$ By setting,  $\tilde z_k^*:=\sum\limits_{j=1}^q\frac{\lambda_j^k}{\lambda_0^k}\big(a_j(u_j^k),b_j(u_j^k)\big)$ for sufficiently large $k\in\N$, we see that $\tilde z_k^*\in C$ for such $k\in \N$ and $\tilde z^*_k=\frac{x^*_k}{\lambda_0^k}-(0_n,1)\to 0-(0_n,1)$ as $k\to\infty.$ In addition, since $C$ is closed, it follows that $-(0_n,1)\in C.$ This means that there exist  $\bar\mu_j\ge 0$ and $\bar u_j\in U_j, j=1,\ldots,p$ such that
\begin{align}\label{1.1}0_n=\sum\limits_{j=1}^q\bar \mu_j a_j(\bar u_j),\quad
-1=\sum\limits_{j=1}^q\bar\mu_j b_j(\bar u_j).\end{align} Besides, by (i), $a_j(u_j)^\top x\le b_j(u_j)$  for all $u_j\in U_j$ and $j=1,\ldots,q$, it follows that \begin{align*}\sum\limits_{j=1}^q\bar \mu_j a_j(\bar u_j)^\top x\le \sum\limits_{j=1}^q\bar\mu_jb_j(\bar u_j).\end{align*} This together with \eqref{1.1} entails that $0\le -1,$ which is absurd, and hence the sequence $\lambda_0^k, k\in \N$ must be bounded. Then, we may assume without loss of generality that $\lambda_0^k\to \lambda_0\ge 0$ as $k\to\infty.$ Similarly, by setting $z_k^*:=\sum\limits_{j=1}^q\lambda_j^k\big(a_j(u^{k}_j),b_j(u^{k}_j)\big)$, we obtain that $z_k^*\in C$ for all $k\in\N$ and $z_k^*\to x^*-\lambda_0(0_n,1)\in C$ as $k\to\infty.$ Therefore, we find $\mu_j\ge 0$ and $u_j\in U_j, j=1,\ldots,p$ such that $$x^*-\lambda_0(0_n,1)=\sum\limits_{j=1}^q\mu_j\big(a_j(u_j),b_j(u_j)\big).$$ This shows that $x^*\in\tilde C,$ and consequently, $\tilde C$ is closed.

Now, invoking Lemma~\ref{Ge-Farkas}, we conclude that $$(-\wp,-r)\in {\rm cl}\tilde C=\tilde C.$$ Then, there exist  $\lambda_0\ge 0$, $\mu_j\ge 0,$ and $u_j:=(u^1_j,\ldots,u^s_j)\in U_j, j=1,\ldots,q$ such that
\begin{align*}-\wp=\sum\limits_{j=1}^q\mu_j\big(a^0_j+\sum^{s}_{i=1}u^i_ja^i_j\big),\quad
-r=\lambda_0+\sum\limits_{j=1}^q\mu_j\big(b^0_j+\sum^{s}_{i=1}u^i_jb^i_j\big).\end{align*}
Using the variable transformations, $\lambda^0_j:=\mu_j\ge 0$ and $\lambda_j^i:=\mu_ju_j^i\in\R, j=1,\ldots,q, i=1,\ldots, s$, we see that $$\wp+\sum\limits_{j=1}^q(\lambda_j^0a^0_j+\sum^{s}_{i=1}\lambda^i_ja^i_j)=0,\quad r+\lambda_0+ \sum\limits_{j=1}^q(\lambda_j^0b^0_j+\sum^{s}_{i=1}\lambda^i_jb^i_j)=0.$$ The later equality means that $-r-\sum\limits_{j=1}^q(\lambda_j^0b^0_j+\sum^{s}_{i=1}\lambda^i_jb^i_j)=\lambda_0\ge 0.$

Let $j\in \{1,\ldots,q\}$ be arbitrary. The relation $u_j\in U_j$ ensures that $A^0_j+\sum_{i=1}^su^i_jA^i_j\succeq 0.$ We  will verify that \begin{align}\label{1.2}\lambda_j^0A^0_j+\sum_{i=1}^s\lambda^i_jA^i_j\succeq 0.\end{align} Indeed, if $\lambda_j^0= 0$, then $\lambda_j^i=0$ for all $i=1,\ldots,s,$  and hence, \eqref{1.2} holds trivially. If $\lambda_j^0\neq 0,$ then $$\lambda_j^0A^0_j+\sum_{i=1}^s\lambda^i_jA^i_j=\lambda_j^0\left(A^0_j+\sum_{i=1}^s\frac{\lambda^i_j}{\lambda_j^0}A^i_j\right)=
\lambda_j^0\left(A^0_j+\sum_{i=1}^su^i_jA^i_j\right)\succeq 0,$$ showing \eqref{1.2} holds, too. Consequently, we obtain (ii).

[{\bf (ii) $\Longrightarrow$ (i)}]  Assume that (ii) holds. It means that there exist  $\lambda_j^0\ge 0, \lambda_j^i\in\R, j=1,\ldots,q, i=1,\ldots,s$ such that $\wp+\sum\limits_{j=1}^q(\lambda_j^0a^0_j+\sum^{s}_{i=1}\lambda^i_ja^i_j)=0, -r-\sum\limits_{j=1}^q(\lambda_j^0b^0_j+\sum^{s}_{i=1}\lambda^i_jb^i_j)\ge 0 \mbox{ and } \lambda_j^0A^0_j+\sum_{i=1}^s\lambda^i_jA^i_j\succeq 0,  j=1,\ldots,q.$ By letting $\lambda_0:=-r-\sum\limits_{j=1}^q(\lambda_j^0b^0_j+\sum^{s}_{i=1}\lambda^i_jb^i_j)$, we obtain that $\lambda_0\ge 0$ and that \begin{align}\label{1.4} -\wp=\sum\limits_{j=1}^q(\lambda_j^0a^0_j+\sum^{s}_{i=1}\lambda^i_ja^i_j),\quad -r=\lambda_0+\sum\limits_{j=1}^q(\lambda_j^0b^0_j+\sum^{s}_{i=1}\lambda^i_jb^i_j),\end{align} and \begin{align}\label{1.3} \lambda_j^0A^0_j+\sum_{i=1}^s\lambda^i_jA^i_j\succeq 0,\quad  j=1,\ldots,q.\end{align}
\medskip

Let  $j\in\{1,\ldots,q\}$ be arbitrary. We claim by \eqref{1.3} that if $\lambda_j^0=0$, then $\lambda_j^i=0$ for all $i=1,\ldots,s.$ Suppose, on the contrary, that $\lambda_j^0=0$ but there exists $i_0\in\{1,\ldots,s\}$ with $\lambda_j^{i_0}\neq 0.$ In this case, \eqref{1.3} becomes $\sum^{s}_{i=1}\lambda^i_jA^i_j\succeq 0$. Let $\bar u_j:= (\bar u^1_j,\ldots,\bar u^s_j)\in U_j.$ It follows by definition that $A^0_j+\sum_{i=1}^s\bar u^i_jA^i_j\succeq 0$ and  $$A^0_j+\sum_{i=1}^s(\bar u^i_j+t\lambda_j^i)A^i_j=\left(A^0_j+\sum_{i=1}^s\bar u^i_jA^i_j\right)+t\sum_{i=1}^s\lambda_j^iA^i_j\succeq 0 \mbox{ for all } t>0.$$ This guarantees that $\bar u_j+t(\lambda_j^1,\ldots,\lambda_j^s)\in U_j$ for all $t>0$, which contradicts the fact that $U_j$ is a compact set. So, our claim must be true.

Next,  take $\hat u_j:=(\hat u_j^1,\ldots,\hat u_j^s)\in U_j$ and define $\tilde u_j:=(\tilde u_j^1,\ldots,\tilde u_j^s)$ with $$\tilde u^i_j:=\begin{cases} \hat u_j^i &\mbox{ if } \lambda_j^0=0,\\ \frac{\lambda_j^i}{\lambda_j^0} &\mbox{ if } \lambda_j^0\neq 0.\end{cases}$$ Then, it can be checked that $$A^0_j+\sum_{i=1}^s\tilde u^i_jA^i_j =\begin{cases} A^0_j+\sum_{i=1}^s\hat u^i_jA^i_j &\mbox{ if } \lambda_j^0=0,\\ \frac{1}{\lambda_j^0}(\lambda_j^0 A^0_j+\sum_{i=1}^s\lambda_j^i A^i_j) &\mbox{ if } \lambda_j^0\neq 0,\end{cases}$$ which shows that $A^0_j+\sum_{i=1}^s\tilde u^i_jA^i_j\succeq 0$, and so, $\tilde u_j\in U_j.$ We now deduce from \eqref{1.4} that  \begin{align*}\label{} (-\wp,-r)=&\lambda_0 (0_n,1)+\sum\limits_{j=1}^q\lambda_j^0\big(a^0_j+\sum^{s}_{i=1}\tilde u^i_ja^i_j,b_j^0+\sum^{s}_{i=1}\tilde u^i_jb^i_j\big)\\=&\lambda_0 (0_n,1)+\sum\limits_{j=1}^q\lambda_j^0\big(a_j(\tilde u_j),b_j(\tilde u_j)\big),\end{align*} which shows that \begin{align}\label{1.5}(-\wp,-r)\in\tilde C,\end{align} where $\tilde C$ is defined as in \eqref{coneC-tilde}.

To prove  (i), let $x\in\R^n$ be such that $a_j(u_j)^\top x\le b_j(u_j)$  for all $u_j\in U_j, j=1,\ldots,q.$ Invoking Lemma~\ref{Ge-Farkas} again, we conclude by \eqref{1.5} that $-\wp^\top x\le -r$ or equivalently, $\wp^\top x\ge r,$ which completes the proof of the theorem. $\hfill\Box$

\medskip In the setting of  $b_j:=0, j=1,\ldots,q$, and $r:=0$, the above result collapses into the so-called {\it semi-infinite Farkas's lemma} given in \cite[Theorem~2.1]{jll-submitted}. The following example shows that  Theorem~\ref{Theo1} may fail without the assumption that the set $C$ is closed. Here, we consider the case of $q:=1$ for the purpose of simplicity.

\begin{Example}\label{} {\bf (The importance of the closed cone regularity)} {\rm Let $U:=\big\{u:=(u^1,u^2)\in\R^2\mid A^0+\sum_{i=1}^2u^iA^i\succeq 0\big\},$ where $$A^0:= \left(
\begin{array}{ccc}
1 & 0&0 \\
0& 1 &0\\ 0& 0 &1\\
\end{array}
\right),\quad A^1:= \left(
\begin{array}{ccc}
0 & 0&1 \\
0& 0 &0\\ 1& 0 &0\\
\end{array}
\right),\quad A^2:= \left(
\begin{array}{ccc}
0 & 0&0 \\
0& 0 &1\\ 0& 1&0\\
\end{array}
\right)\cdot$$   Let $a:\R^2\to \R^2$ and $b:\R^2\to\R$ be affine mappings defined respectively by $a(u):=a^0+\sum^{2}_{i=1}u^ia^i$ and $b(u):=b^0+\sum^{2}_{i=1}u^ib^i$ for $u:=(u^1, u^2)\in \R^2$ with $a^0:=(0,1), a^1:=(1,0), a^2:=(0,1)\in \R^2, b^0=b^1=b^2:=0\in\R.$

In this setting, it is easy to see that $U=\{u:=(u^1,u^2)\in\R^2\mid (u^1)^2+(u^2)^2\le 1\}.$ Let $x:=(x^1,x^2)\in\R^2$ be such that $a(u)^\top x\le b(u)$  for all $u:=(u^1,u^2)\in U.$ It means that $x$ satisfies \begin{align}\label{2.6}u^1x^1+(u^2+1)x^2\le 0,\quad \forall u:=(u^1,u^2)\in U.\end{align} Now, choose $\bar x:=(0,-1), \wp:=(1,0)$ and $r:=0.$ We see that $\bar x$ satisfies \eqref{2.6} and $\wp^\top\bar x\ge r.$ It means that we have the condition (i) of Theorem~\ref{Theo1}.
\medskip

However, condition (ii) of Theorem~\ref{Theo1} fails. Indeed, assume on the contrary that there exist $\lambda^0\ge 0$ and $\lambda^1, \lambda^2\in\R$ such that \begin{align*}\wp+\lambda^0a^0+\lambda^1a^1+\lambda^2a^2=(0,0),\; -r-(\lambda^0b^0+\lambda^1b^1+\lambda^2b^2)\ge 0,\; \lambda^0A^0+\lambda^1A^1+\lambda^2A^2\succeq 0.\end{align*} It reduces to the following expression \begin{align*}\lambda^1=-1,\; \lambda^2=-\lambda^0, \; (\lambda^0)^2\ge (\lambda^1)^2+(\lambda^2)^2,\end{align*} which is absurd.

Consequently, the conclusion  of Theorem~\ref{Theo1} fails to hold. The reason is that  the  cone $C:={\rm cone}\big\{(u^1,u^2+1,0)\mid (u^1)^2+(u^2)^2\le 1\big\}$ is not closed. To see this, just  take the sequence $z_k:=(1,\frac{1}{k},0)=k(\frac{1}{k},\frac{1}{k^2},0)\in C$ for  $k\in\N$. It is clear that $z_k\to  z_0:=(1,0,0)$ as $k\to\infty,$ and $z_0\notin C.$
} \end{Example}

We now   provide some sufficient criteria which guarantee that the  cone $C$ in Theorem~\ref{Theo1} is closed.

\begin{Proposition}\label{Pro1}{\bf (Sufficiency for the closed cone regularity condition)} Assume that one of the following conditions holds:

{\rm (i)\ } $U_j, j=1,\ldots,q$ are  bounded polyhedra $(\mbox{i.e.,   polytopes});$

{\rm (ii)\ } $U_j, j=1,\ldots,q$ are  compact  and the Slater condition holds, i.e., there is $\hat x\in \R^n$ such that  \begin{align}\label{Slater} a_j(u_j)^\top \hat x< b_j(u_j), \;\forall u_j\in U_j,\; j=1,\ldots,q.\end{align}

Then, the  cone  $C:={\rm cone}\big\{\big(a_j(u_j),b_j(u_j)\big)\mid u_j\in U_j, j=1,\ldots,q\big\}$  is closed.
\end{Proposition}
{\bf Proof.} For each $j\in\{1,\ldots,q\},$ define an affine map $F_j:\R^s\to \R^n\times\R$ by $$F_j(u_j):=[a_j(u_j),b_j(u_j)\big].$$ Then, $F_j(U_j):=\bigcup\limits_{u_j\in U_j}F(u_j)=\bigcup\limits_{u_j\in U_j}[a_j(u_j),b_j(u_j)\big]$, and so it follows that $$C={\rm cone}\left(\bigcup\limits_{j=1}^{q}F_j(U_j)\right).$$

(i) Since $U_j$ is a polytope for $j\in\{1,\ldots,q\},$  the set $F_j(U_j)$ is a polytope as well (see e.g., \cite[Proposition~1.29]{Bruns-Gu09}), and so $F_j(U_j)$ is the convex hull of some finite set (see e.g., \cite[Theorem~1.26]{Bruns-Gu09}). It means that for each $j\in\{1,\ldots,q\},$ there exists $k_j\in\N$ and $(a_j^i,b^i_j)\in\R^n\times\R, i=1,\ldots,k_j$ such that $F_j(U_j)={\rm conv}\{(a_j^i,b^i_j)\mid i=1,\ldots,k_j\}.$ It follows that ${\rm conv}\left(\bigcup\limits_{j=1}^{q}F_j(U_j)\right)={\rm conv}\{(a_j^i,b^i_j)\mid i=1,\ldots,k_j,\; j=1,\ldots,q\}$ is a polytope. So, $C={\rm cone}\left(\bigcup\limits_{j=1}^{q}F_j(U_j)\right)={\rm cone}\left[{\rm conv}\left(\bigcup\limits_{j=1}^{q}F_j(U_j)\right)\right]$ is a polyhedral cone, and hence, is closed (see e.g., \cite[Proposition~3.9]{Mor-Nam-14}).

(ii) Since $U_j$ is compact  and  $F_j$ is  affine (and thus, continuous) for $j\in\{1,\ldots,q\}$, the set $F_j(U_j)$ is a compact set. It entails that $\left(\bigcup\limits_{j=1}^{q}F_j(U_j)\right)$ is a compact set as well. Note further that  $F_j(U_j)$ is a convex set for each $j\in\{1,\ldots,q\}$ due to the convexity of $U_j$ (see e.g., \cite[Proposition~1.23]{Mor-Nam-14}).
Assume now that the Slater condition~\eqref{Slater} holds. We claim that \begin{align}\label{1.8}0_{n+1}\notin {\rm conv}\left(\bigcup\limits_{j=1}^{q}F_j(U_j)\right).\end{align} If this is not the case, then there exist $\bar u_j\in U_j, \bar\mu_j\ge 0, j=1,\ldots,q$ with $\sum\limits_{j=1}^q\bar\mu_j=1$ such that \begin{align}\label{1.7}0_{n+1}=\sum\limits_{j=1}^q\bar\mu_j\big[a_j(\bar u_j),b_j(\bar u_j)\big].\end{align} Moreover, we get by \eqref{Slater} that $$\sum\limits_{j=1}^q\bar\mu_ja_j(\bar u_j)^\top \hat x<\sum\limits_{j=1}^q\bar\mu_jb_j(\bar u_j).$$ It together with \eqref{1.7} gives a contradiction, and so, \eqref{1.8} holds. According to  \cite[Proposition~1.4.7]{Hiriart-Urruty93}, we conclude that the cone $C={\rm cone}\left(\bigcup\limits_{j=1}^{q}F_j(U_j)\right)$ is closed. $\hfill\Box.$

\medskip
As an application of Theorem 2.2, we obtain a complete dual characterization of robust optimality of the uncertain lower level linear program of  the  bilevel polynomial problem.

\medskip
Let $x\in\R^m$ and consider the uncertain lower level linear program~\eqref{LNP}, given by \begin{align}\label{LNP}\min_{z\in\R^n}{\{c_0^\top x+d_0^\top z\mid c_j^\top x +a_j(u_j)^\top z\le b_j(u_j),\; j=1,\ldots,q\}},\tag{LNP}\end{align} where $u_j\in U_j, \; j=1,\ldots,q$ with the parameter sets $U_j:=\big\{(u^1_j,\ldots, u^s_j)\in\R^s\mid A^0_j+\sum_{i=1}^su^i_jA^i_j\succeq 0\big\},\; j=1,\ldots, q,$  as in \eqref{U-set}, and the affine mappings $a_j:\R^s\to \R^n, b_j:\R^s\to\R, j=1,\ldots,q$ are  given respectively by $a_j(u_j):=a^0_j+\sum^{s}_{i=1}u^i_ja^i_j, b_j(u_j)=b^0_j+\sum^{s}_{i=1}u^i_jb^i_j$ as in \eqref{affine-maps}, as well as  $c_0\in\R^m, d_0\in\R^n, c_j \in \R^m,    j=1,\ldots,q$ fixed.

\medskip
Following the robust optimization approach \cite{robustbook,goberna-jeya-siam,goberna-jeya-MP,jeya-optimL}, the {\it robust} counterpart of the problem~\eqref{LNP} is defined by
 \begin{align}\label{RLNP}\min_{z\in\R^n}{\{c_0^\top x+d_0^\top z\mid c_j^\top x +a_j(u_j)^\top z\le b_j(u_j),\;\forall u_j\in U_j,\; j=1,\ldots,q\}}.\tag{RLNP}\end{align}

 As usual, we say that $y\in \R^n$ is a {\it robust solution} of the problem~\eqref{LNP} (see e.g., \cite{robustbook}) if it is an optimal solution of the problem~\eqref{RLNP}, i.e., $y\in Y(x):={\rm argmin}_{z\in\R^n}\{c_0^\top x+d_0^\top z\mid c_j^\top x +a_j(u_j)^\top z\le b_j(u_j),\; \forall u_j\in U_j,\;  j=1,\ldots,q\}.$

 \medskip
 The next result establishes a characterization for robust solutions of the uncertain linear optimization problem~\eqref{LNP}.

\begin{Theorem}\label{Theo2-LN}{\bf (Characterization for robust solutions of problem~\eqref{LNP})} Let $x\in\R^m$, and let the cone $C(x):={\rm cone}\big\{\big(a_j(u_j),b_j(u_j)-c_j^\top x\big)\mid u_j\in U_j, j=1,\ldots,q\big\}$ be closed. Then,  $y\in Y(x)$ if and only if \begin{align*}({\rm I})\begin{cases}&y\in\R^n, c_j^\top x +a_j(u_j)^\top y\le b_j(u_j),\;\forall u_j\in U_j,\; j=1,\ldots,q \mbox{ and }\\ &\exists  \lambda_j^0\ge 0, \lambda_j^i\in\R, j=1,\ldots,q, i=1,\ldots,s \mbox{ such that }\\&d_0+\sum\limits_{j=1}^q(\lambda_j^0a^0_j+\sum^{s}_{i=1}\lambda^i_ja^i_j)=0, \; -d_0^\top y-\sum\limits_{j=1}^q(\lambda_j^0b^0_j-\lambda_j^0c_j^\top x+\sum^{s}_{i=1}\lambda^i_jb^i_j)\ge 0 \\&\mbox{and } \lambda_j^0A^0_j+\sum_{i=1}^s\lambda^i_jA^i_j\succeq 0,  j=1,\ldots,q.\end{cases}\end{align*}
\end{Theorem}
{\bf Proof.} Let $y\in Y(x).$ This means that \begin{align*}\label{} y\in\R^n, \; c_j^\top x +a_j(u_j)^\top y\le b_j(u_j),\;\forall u_j\in U_j,\; j=1,\ldots,q\end{align*} and \begin{align}\label{1.6-sua}d_0^\top y\le d_0^\top z\mbox{ for all } z\in\R^n \mbox{ satisfying } c_j^\top x+ a_j(u_j)^\top z\le b_j(u_j),\;\forall u_j\in U_j,\; j=1,\ldots,q.\end{align} Letting $r:=d_0^\top y,$  \eqref{1.6-sua} amounts to the assertion that for each  \begin{align*}&z\in\R^n,\; (a^0_j+\sum^{s}_{i=1}u^i_ja^i_j)^\top z\le (b^0_j-c_j^\top x)+\sum^{s}_{i=1}u^i_jb^i_j,\;  \forall u_j:=(u^1_j,\ldots, u^s_j)\in U_j, j=1,\ldots,q\\&\Longrightarrow d_0^\top z\ge r.\end{align*} Since the cone $C(x)$ is closed, applying  Theorem~\ref{Theo1}, we find $ \lambda_j^0\ge 0, \lambda_j^i\in\R, j=1,\ldots,q, i=1,\ldots,s $  such that \begin{align*}&d_0+\sum\limits_{j=1}^q(\lambda_j^0a^0_j+\sum^{s}_{i=1}\lambda^i_ja^i_j)=0, \; -r-\sum\limits_{j=1}^q(\lambda_j^0b^0_j-\lambda_j^0c_j^\top x+\sum^{s}_{i=1}\lambda^i_jb^i_j)\ge 0 \\&\mbox{and } \lambda_j^0A^0_j+\sum_{i=1}^s\lambda^i_jA^i_j\succeq 0,  j=1,\ldots,q.\end{align*} So, we obtain (I).

Conversely, assume that (I) holds. Then, there exists $y\in\R^n$ such that $$c_j^\top x +a_j(u_j)^\top y\le b_j(u_j),\;\forall u_j\in U_j,\; j=1,\ldots,q.$$ This means that $y\in\R^n$ is a feasible point of problem~\eqref{RLNP}.
Let $z\in\R^n$ be such that $c_j^\top x+a_j(u_j)^\top z\le b_j(u_j),\;\forall u_j\in U_j,\; j=1,\ldots,q$ or equivalently, $$ (a^0_j+\sum^{s}_{i=1}u^i_ja^i_j)^\top z\le (b^0_j-c_j^\top x)+\sum^{s}_{i=1}u^i_jb^i_j,\;  \forall u_j:=(u^1_j,\ldots, u^s_j)\in U_j, j=1,\ldots,q.$$ By (I), there exist $\lambda_j^0\ge 0, \lambda_j^i\in\R, j=1,\ldots,q, i=1,\ldots,s \mbox{ such that }$ \begin{align*}&d_0+\sum\limits_{j=1}^q(\lambda_j^0a^0_j+\sum^{s}_{i=1}\lambda^i_ja^i_j)=0, \; -d_0^\top y-\sum\limits_{j=1}^q(\lambda_j^0b^0_j-\lambda_j^0c_j^\top x+\sum^{s}_{i=1}\lambda^i_jb^i_j)\ge 0 \\&\mbox{and } \lambda_j^0A^0_j+\sum_{i=1}^s\lambda^i_jA^i_j\succeq 0,  j=1,\ldots,q.\end{align*} Due to the fact that  $C(x)$ is closed, invoking Theorem~\ref{Theo1} again, we conclude that $$d_0^\top z\ge d_0^\top y.$$ So, $y\in Y(x),$ which finishes the proof of the theorem. $\hfill\Box$

\section{Uncertain Bilevel Problems \& Box Data Uncertainty}
\setcounter{equation}{0}
This section is devoted to examining a  bilevel  polynomial optimization problem with uncertain linear constraints. Let us first recall that a real-valued function $f:\R^n\to \R$ is {\it coercive} on $\R^n$ if $\liminf\limits_{||x||\to\infty}{f(x)}=+\infty$. In particular,  a convex polynomial $f$  is coercive on $\R^n$ if there exists $\bar x\in\R^n$ such that the Hessian $\nabla^2f(\bar x)$ is positive definite (see e.g., \cite[Lemma~3.1]{jeya-Son-Li-ORL14}). Numerically checkable sufficient conditions for the coercivity of nonconvex polynomials have also been given in \cite{JLL}.

\medskip
Denote by $\R[x]$ the ring of polynomials in $x$ with real coefficients.  One says that (cf.~\cite{Lasserre-book-09})  $f\in\R[x]$ is {\it sums-of-squares} (SOS) if there exist polynomials $f_j\in\R[x], j=1,\ldots, r$  such that $f=\sum_{j=1}^rf_j^2.$ The set of all (SOS) polynomials in $x$ is denoted by $\Sigma^2[x].$ Given polynomials $\{g_1,\ldots,g_r\}\subset \R[x],$ the notation $\mathbf{M}(g_1,\ldots,g_r)$ stands for the set of polynomials generated by $\{g_1,\ldots,g_r\},$ i.e., \begin{align}\label{Set-polynomial} \mathbf{M}(g_1,\ldots,g_r):=\{\sigma_0+\sigma_1g_1+\ldots+\sigma_rg_r\mid \sigma_j\in \Sigma^2[x], j=0,1,\ldots,r\}.\end{align}

The set (cf.~\cite{Demel-Nie07}) $\mathbf{M}(g_1,\ldots,g_r)$ is {\it archimedean} if there exists $h\in \mathbf{M}(g_1,\ldots,g_r)$ such that the set $\{x\in\R^n\mid h(x)\ge 0\}$ is compact.

\medskip
The following lemma of Putinar (cf.~\cite{Lasserre-book-09,Putinar93}), which
provides a positivity representation for a polynomial over a system of
polynomial inequalities under the archimedean property, can be viewed as a
polynomial analog of Farkas's Lemma.

\begin{Lemma}\label{Putinar} {\rm (Putinar's Positivstellensatz \cite{Putinar93})} Let $f, g_j\in\R[x], j=1,\ldots, r.$ Suppose that $\mathbf{M}(g_1,\ldots,g_r)$ is archimedean. If $f(x)>0$ for all $x\in K:=\{y\in\R^n\mid g_j(y)\ge 0, j=1,\ldots,r\}$, then $f\in \mathbf{M}(g_1,\ldots,g_r),$ i.e., there exist $\sigma_j\in \Sigma^2[x], j=0,1,\ldots,r$ such that $f=\sigma_0+\sum_{j=1}^r\sigma_jg_j.$
\end{Lemma}

Let $f:\R^m\times\R^n\rightarrow\R$ be a real polynomial. We consider  the {\it  bilevel polynomial}  optimization problem with {\it uncertain linear} constraints as \begin{eqnarray*}
 {\rm (P)}&\displaystyle\min_{(x,y)\in \R^{m}\times \R^{n}}&
f(x,y) \\
&\mbox{ subject to }&  \tilde a_i^\top x+\tilde b_i^\top y\le \tilde c_i,\; i=1,\ldots,l, \\
&&y\in Y(x, \hat a_1,\hat b_1, \ldots, \hat a_q,\hat b_q),
\end{eqnarray*} where $(\tilde a_i, \tilde b_i, \tilde c_i)\in \widetilde U_i, i=1,\ldots,l$ and $(\hat a_j,\hat b_j)\in \widehat U_j, j=1,\ldots, q$  are {\it uncertain} and  $Y(x, \hat a_1,\hat b_1, \ldots, \hat a_q,\hat b_q):={\rm argmin}_{z\in\R^n}\{c_0^\top x+d_0^\top z\mid c_j^\top x +\hat a_j^\top z\le \hat b_j,\;  j=1,\ldots,q\}$ denotes the optimal solution set of the uncertain lower-level optimization problem \begin{align}\label{LP-origi}\min_{z\in\R^n}{\{c_0^\top x+d_0^\top z\mid c_j^\top x +\hat a_j^\top z\le \hat b_j,\; j=1,\ldots,q\}}.\end{align}
In the above data,  $c_0\in\R^m, d_0\in\R^n, c_j \in \R^m,    j=1,\ldots,q$ are fixed, the  {\it uncertainty} sets $\widetilde U_i, i=1,\ldots,l$ are boxes given by  \begin{align*}\label{}\widetilde U_i:=[\underline{a}_i,\overline{a}_i]\times [\underline{b}_i,\overline{b}_i]\times[\underline{c}_i,\overline{c}_i], \; i=1,\ldots,l\end{align*} with $\underline{a}_i:=(\underline{a}_i^1,\ldots,\underline{a}_i^m), \overline{a}_i:=(\overline{a}_i^1,\ldots,\overline{a}_i^m)\in\R^m, \underline{a}_i^k\le \overline{a}_i^k, k=1,\ldots,m, \underline{b}_i:=(\underline{b}_i^1,\ldots,\underline{b}_i^n), \overline{b}_i:=(\overline{b}_i^1,\ldots,\overline{b}_i^n)\in\R^n, \underline{b}_i^k\le \overline{b}_i^k, k=1,\ldots,n,$ and $\underline{c}_i, \overline{c}_i\in\R, \underline{c}_i\le \overline{c}_i$ for $i=1,\ldots,l,$ while the  {\it uncertainty} sets $\widehat U_j, j=1,\ldots, q$ are given by $$\widehat U_j:=\Big\{(a^0_j+\sum^{s}_{i=1}u^i_ja^i_j,b^0_j+\sum^{s}_{i=1}u^i_jb^i_j)\mid  u_j:=(u^1_j,\ldots, u^s_j)\in U_j\Big\},\; j=1,\ldots,q$$ with $U_j:=V_s:=[-\gamma_1,\gamma_1]\times\cdots\times[-\gamma_s,\gamma_s], \gamma_i>0, i=1,\ldots, s$ and $a_j^i\in \R^n, b_j^i\in\R, i=0,1,\ldots,s, j=1,\ldots,q$ fixed.

\medskip
As shown in the previous sections, by considering the affine mappings $a_j:\R^s\to \R^n, b_j:\R^s\to\R, j=1,\ldots,q$  given respectively by $a_j(u_j):=a^0_j+\sum^{s}_{i=1}u^i_ja^i_j, b_j(u_j)=b^0_j+\sum^{s}_{i=1}u^i_jb^i_j$ as in \eqref{affine-maps}, the uncertain lower-level optimization problem~\eqref{LP-origi} can be formulated equivalently as \begin{align}\label{LP}\min_{z\in\R^n}{\{c_0^\top x+d_0^\top z\mid c_j^\top x +a_j(u_j)^\top z\le b_j(u_j),\; j=1,\ldots,q\}},\tag{LP}\end{align} where $u_j\in U_j, j=1,\ldots,q.$ In the formulation of~\eqref{LP}, the parameter sets $U_j, j=1,\ldots, q$ play the role of uncertainty sets.

\medskip
Now, the {\it robust} counterpart of the problem~(P) is defined by
\begin{align}\label{RP}  \min_{(x,y)\in \R^m\times\R^n}{\big\{f(x,y)} \mid  y\in Y(x),\;\tilde a_i^\top x+\tilde b_i^\top y\le \tilde c_i,\; \forall (\tilde a_i, \tilde b_i, \tilde c_i)\in\widetilde U_i,\; i=1,\ldots,l\big\}, \tag{RP}\end{align}
where  $Y(x):={\rm argmin}_{z\in\R^n}\{c_0^\top x+d_0^\top z\mid c_j^\top x +a_j(u_j)^\top z\le b_j(u_j),\; \forall u_j\in U_j,\;  j=1,\ldots,q\}.$ Note that in the robust counterpart \eqref{RP} the
uncertain constraint inequalities of both the lower-level and upper-level problems  are enforced for every possible value of the data within
the uncertainty sets $U_j, j=1,\ldots,q$ and $\widetilde U_i,\; i=1,\ldots,l.$

\medskip
Given $x\in\R^m$, let us first focus on  the robust counterpart of the lower-level optimization problem~\eqref{LP} given by \begin{align}\label{RLP}\min_{z\in\R^n}{\{c_0^\top x+d_0^\top z\mid c_j^\top x +a_j(u_j)^\top z\le b_j(u_j),\;\forall u_j\in U_j,\; j=1,\ldots,q\}}.\tag{RLP}\end{align}

 Similar to the definition of robust solutions for the uncertain linear optimization problem~\eqref{LNP}, a point $y\in \R^n$ is a  robust solution of the  lower-level optimization problem~\eqref{LP} if it is an optimal solution of the problem~\eqref{RLP}, i.e., $y\in Y(x):={\rm argmin}_{z\in\R^n}\{c_0^\top x+d_0^\top z\mid c_j^\top x +a_j(u_j)^\top z\le b_j(u_j),\; \forall u_j\in U_j,\;  j=1,\ldots,q\}.$

\medskip
We observe here that for each $ j\in\{1,\ldots, q\},$ put \begin{align}\label{2.8} A^0_j:= \left(
\begin{array}{cc}
E_0 & 0 \\
0& E_0 \\
\end{array}
\right),\quad A^i_j:= \left(
\begin{array}{cc}
E_i & 0 \\
0& -E_i \\
\end{array}
\right),\; i=1,\ldots, s,\end{align} where $E_0$ is the $(s\times s)$ diagonal matrix with the diagonal entries, say $\gamma_i>0, i=1,\ldots, s,$ and $E_i$ is the $(s\times s)$ diagonal matrix with one in the $(i,i)$th entry and zeros elsewhere. Then,  $$A^0_j+\sum_{i=1}^su^i_jA^i_j=\begin{bmatrix}
    \gamma_1+u^1_j&0&\cdots&0&0&\cdots & 0 \\
    0&\gamma_2+u^2_j&\cdots&0&0&\cdots & 0\\
    \vdots&\vdots&\vdots&\vdots&\vdots&\vdots & \vdots\\
    0&0&\gamma_s+u^s_j&0&0&\cdots & 0\\
    0&0&0&\gamma_1-u^1_j&0&\cdots & 0\\
 \vdots&\vdots&\vdots&\vdots&\vdots&\vdots & \vdots\\
  0&0&\cdots&0&0&\cdots & \gamma_s-u^s_j
  \end{bmatrix},$$ and therefore, we have
\begin{align*}\label{}&\{u_j:=(u^1_j,\ldots, u^s_j)\in\R^s\mid A^0_j+\sum_{i=1}^su^i_jA^i_j\succeq 0\big\}\notag\\&=\big\{u_j:=(u^1_j,\ldots, u^s_j)\in\R^s\mid
\gamma_i+u_j^i\ge 0, \gamma_i-u_j^i\ge 0, i=1,\ldots, s\big\}\\&=\big\{u_j:=(u^1_j,\ldots, u^s_j)\in\R^s\mid |u^i_j|\le \gamma_i, i=1,\ldots, s\big\}\notag\\&=V_s,\; j=1,\ldots, q,\notag\end{align*} which shows how the box $V_s$ can be expressed in terms of  spectrahedra in \eqref{U-set}. From now on, we denote by $\{\check{u}_k:=(\check{u}^{1}_k,\ldots,\check{u}^{s}_k)\in\R^s\mid k=1,\ldots, 2^s\}$ the {\it extreme points} of the box $V_s$.

\medskip
The following corollary establishes a characterization for robust solutions of the lower-level optimization problem~\eqref{LP}.

\begin{Corollary}\label{Theo2} {\bf (Tractable characterization for robust solutions of~\eqref{LP})} Let $x\in\R^m$. Then, $y\in Y(x)$ if and only if  \begin{align*}({\rm I})\begin{cases}&y\in\R^n, c^\top_jx +(a_j^0+\sum^{s}_{i=1}\check{u}^{i}_ka^i_j)^\top y -(b_j^0+ \sum^{s}_{i=1}\check{u}^{i}_kb^i_j)\le 0, k=1,\ldots, 2^s, j=1,\ldots,q, \\ &\exists \mu_0> 0, \mu_j\ge 0, \mu_j^i\in\R, j=1,\ldots,q, i=1,\ldots,s \mbox{ such that } \sum\limits_{j=0}^q(\mu_j)^2+\sum\limits_{j=1}^q\sum\limits_{i=1}^s(\mu_j^i)^2=1, \\& \mu_0d_0+\sum\limits_{j=1}^q(\mu_ja^0_j+\sum^{s}_{i=1}\mu^i_ja^i_j)=0, \; -\mu_0d_0^\top y-\sum\limits_{j=1}^q(\mu_jb^0_j-\mu_jc_j^\top x+\sum^{s}_{i=1}\mu^i_jb^i_j)\ge 0 \\&\mbox{and } (\mu_j\gamma_i)^2-(\mu^i_j)^2\ge 0, i=1,\ldots, s,  j=1,\ldots,q.\end{cases}\end{align*}
\end{Corollary}
{\bf Proof.} Since $U_j:=V_s, j=1,\ldots,q$ are polytopes, the  cone  $$C(x):={\rm cone}\big\{\big(a_j(u_j),b_j(u_j)-c_j^\top x\big)\mid u_j\in U_j, j=1,\ldots,q\big\}$$ is closed by virtue of Proposition~\ref{Pro1}(i). According to Theorem~\ref{Theo2-LN}, we conclude that $y\in Y(x)$ if and only if the following conditions hold: \begin{align}&\label{1.13a} y\in\R^n,\; c_j^\top x +a_j(u_j)^\top y\le b_j(u_j),\;\forall u_j\in U_j,\; j=1,\ldots,q \mbox{ and }\\ &\exists  \lambda_j^0\ge 0, \lambda_j^i\in\R, j=1,\ldots,q, i=1,\ldots,s \mbox{ such that }\notag\\&d_0+\sum\limits_{j=1}^q(\lambda_j^0a^0_j+\sum^{s}_{i=1}\lambda^i_ja^i_j)=0, \; -d_0^\top y-\sum\limits_{j=1}^q(\lambda_j^0b^0_j-\lambda_j^0c_j^\top x+\sum^{s}_{i=1}\lambda^i_jb^i_j)\ge 0 \notag\\&\mbox{and } \lambda_j^0A^0_j+\sum_{i=1}^s\lambda^i_jA^i_j\succeq 0,  j=1,\ldots,q.\label{2.7}\end{align}
It is easy to see that \eqref{1.13a} amounts  to the following one $$y\in\R^n, \max_{u_j\in U_j=V_s}{\big\{ \langle (a^{1\top}_jy-b^1_j,\ldots, a_j^{s\top}y-b^s_j),u_j \rangle +a^{0\top}_jy+c^\top_jx-b_j^0 \big\}}\le 0,  j=1,\ldots,q,$$ which in turn  is equivalent to the assertion
$$y\in\R^n, \check{u}^{\top}_k \big(a^{1\top}_jy-b^1_j,\ldots, a_j^{s\top}y-b^s_j\big) +a^{0\top}_jy+c^\top_jx-b_j^0 \le 0, k=1,\ldots, 2^s, j=1,\ldots,q,$$ where $\{\check{u}_k:=(\check{u}^{1}_k,\ldots,\check{u}^{s}_k)\in\R^s\mid k=1,\ldots, 2^s\}$ are the extreme points of the box $V_s$ denoted as above. So, \eqref{1.13a} becomes $$y\in\R^n, c^\top_jx +(a_j^0+\sum^{s}_{i=1}\check{u}^{i}_ka^i_j)^\top y -(b_j^0+ \sum^{s}_{i=1}\check{u}^{i}_kb^i_j)\le 0, k=1,\ldots, 2^s, j=1,\ldots,q.$$ Next, we show that \eqref{2.7} becomes \begin{align}\label{1.15a}(\lambda_j^0\gamma_i)^2-(\lambda^i_j)^2\ge 0, i=1,\ldots, s,  j=1,\ldots,q\end{align} under our setting. Indeed, for each $j\in\{1,\ldots,q\},$  consider the matrices $A^i_j, i=0,1,\ldots,s$ as given in \eqref{2.8}.  Then, we have $$\lambda_j^0A^0_j+\sum_{i=1}^s\lambda^i_jA^i_j=\begin{bmatrix}
   \lambda_j^0\gamma_1+\lambda^1_j&0&\cdots&0&0&\cdots & 0 \\
    0&\lambda_j^0\gamma_2+\lambda^2_j&\cdots&0&0&\cdots & 0\\
    \vdots&\vdots&\vdots&\vdots&\vdots&\vdots & \vdots\\
    0&0&\lambda_j^0\gamma_s+\lambda^s_j&0&0&\cdots & 0\\
    0&0&0&\lambda_j^0\gamma_1-\lambda^1_j&0&\cdots & 0\\
 \vdots&\vdots&\vdots&\vdots&\vdots&\vdots & \vdots\\
  0&0&\cdots&0&0&\cdots & \lambda_j^0\gamma_s-\lambda^s_j
  \end{bmatrix}.$$
This follows  that \begin{align*} \label{} \lambda_j^0A^0_j+\sum_{i=1}^s\lambda^i_jA^i_j\succeq 0 \Leftrightarrow \lambda_j^0\gamma_i+\lambda_j^i\ge 0, \; \lambda_j^0\gamma_i-\lambda_j^i\ge 0,\;  i=1,\ldots, s,  \end{align*} which proves that \eqref{2.7} amounts to \eqref{1.15a}.

So, we come to the assertion that $y\in Y(x)$ if and only if \begin{align*}({\rm II})\begin{cases}&y\in\R^n, c^\top_jx +(a_j^0+\sum^{s}_{i=1}\check{u}^{i}_ka^i_j)^\top y -(b_j^0+ \sum^{s}_{i=1}\check{u}^{i}_kb^i_j)\le 0, k=1,\ldots, 2^s, j=1,\ldots,q, \\ &\exists  \lambda_j^0\ge 0, \lambda_j^i\in\R, j=1,\ldots,q, i=1,\ldots,s \mbox{ such that }\\&d_0+\sum\limits_{j=1}^q(\lambda_j^0a^0_j+\sum^{s}_{i=1}\lambda^i_ja^i_j)=0, \; -d_0^\top y-\sum\limits_{j=1}^q(\lambda_j^0b^0_j-\lambda_j^0c_j^\top x+\sum^{s}_{i=1}\lambda^i_jb^i_j)\ge 0 \\&\mbox{and } (\lambda_j^0\gamma_i)^2-(\lambda^i_j)^2\ge 0, i=1,\ldots, s,  j=1,\ldots,q.\end{cases}\end{align*}

Now,  it is easy to see that (II) entails  (I) with $\mu_0:=\frac{1}{\sqrt{1+\sum\limits_{j=1}^q(\lambda_j^0)^2+\sum\limits_{j=1}^q\sum\limits_{i=1}^s(\lambda_j^i)^2}}$ and $\mu_j:=\frac{\lambda_j^0}{\sqrt{1+\sum\limits_{j=1}^q(\lambda_j^0)^2+\sum\limits_{j=1}^q\sum\limits_{i=1}^s(\lambda_j^i)^2}}, \mu_j^i:=\frac{\lambda_j^i}{\sqrt{1+\sum\limits_{j=1}^q(\lambda_j^0)^2+\sum\limits_{j=1}^q\sum\limits_{i=1}^s(\lambda_j^i)^2}},  j=1,\ldots,q, i=1,\ldots, s$. Conversely, suppose that (I) holds.  By letting $\lambda_j^0:=\frac{\mu_j}{\mu_0}, \lambda_j^i:=\frac{\mu_j^i}{\mu_0}, j=1,\ldots, q, i=1,\ldots, s$, we arrive at (II). So,  (I) and (II) are equivalent, which finishes the proof. $\hfill\Box$

\medskip
 To proceed further, we should define  concepts of robust feasible/solutions for the bilevel polynomial optimization problem~(P) under uncertainty of both levels.

 \begin{Definition}{\rm (i)\ We say that  $(\bar x,\bar y)\in\R^m\times\R^n$ is a {\it  robust feasible} point of problem~(P) if  it satisfies  $$ \bar y\in Y(\bar x),\;\tilde a_i^\top \bar x+\tilde b_i^\top \bar y\le \tilde c_i,\; \forall (\tilde a_i, \tilde b_i, \tilde c_i)\in\widetilde U_i,\; i=1,\ldots,l,$$ or equivalently, it is a feasible point of its  robust counterpart~\eqref{RP}.\\
 (ii)\  Let $(\bar x,\bar y)\in\R^m\times\R^n$ be a  robust feasible point of problem~(P).   We say that $(\bar x,\bar y)$ is a  {\it global  robust solution} of  problem~(P) if $f(\bar x,\bar y)\le f(x,y)$ for every  robust feasible point $(x,y)$ of problem~(P), or equivalently, it is  a global solution of its  robust counterpart~\eqref{RP}.\\
 (iii)\  We say that the  problem~(P) satisfies the  {\it lower-level Slater condition}  (LSC) if  for each $x\in\R^m$ there exists $z\in\R^n$ such that \begin{align}\label{LLSC} c_j^\top x +a_j(u_j)^\top z< b_j(u_j),\;\forall u_j\in U_j,\; j=1,\ldots,q.\end{align}
  }\end{Definition}

  We are now ready to  provide  a characterization for global robust solutions of the bilevel polynomial  optimization problem with uncertain linear constraints~(P). In the following theorem, we will use notation as before, and in addition, we put $d_0:=(d_0^1,\ldots,d_0^n), a_j^i:=(a^{i1}_j,\ldots,a^{in}_j)\in\R^n, i=0,1,\ldots, s, j=1,\ldots, q$, and let $\{(\check{a}^{k}_i,\check{b}^{k}_i,\check{c}^{k}_i)\in\R^m\times\R^n\times\R\mid k=1,\ldots,2^{m+n+1}\}$  denote the  extreme points of the box $\widetilde U_i$ for $i=1,\ldots,l.$

\begin{Theorem}\label{Theo5}{\bf (Characterization for  global robust solutions of~(P))} Let $f$ be coercive on $\R^m\times\R^n$, and let the {\rm (LSC)} in \eqref{LLSC} be satisfied. Let $(\bar x,\bar y)\in\R^m\times\R^n$ be a  robust feasible point of problem~{\rm(P)}, and let $\kappa\in\R$ be such that $\kappa\ge f(\bar x,\bar y).$ Then, $(\bar x,\bar y)$ is a global robust solution of problem~{\rm (P)}  if and only if for any $\epsilon>0$, there exist sums-of-squares polynomials  $\zeta, \sigma_0, \sigma_i\in \Sigma^2[x,y,\mu], i=1,\ldots,L:=l2^{m+n+1}+q(2^s+s+1)+2,$ where $\mu:=(\mu_0,\mu_1,\ldots,\mu_q,\mu_1^1,\ldots,\mu^1_q,\ldots,\mu^s_1,\ldots,\mu^s_q)\in\R^{q(s+1)+1}$ and real polynomials $\xi_j\in\R[x,y,\mu],   j=1,\ldots, n+1$ such that \begin{align}\label{1.28}f-\sum\limits_{i=1}^{L}\sigma_ig_i-\sum\limits_{j=1}^{n+1}\xi_jh_j-\zeta\big(\kappa-f\big)-f(\bar x,\bar y)+\epsilon=\sigma_0\end{align} with \begin{align*}g_i(x,y,\mu):=&- (\check{a}_k^{i-(k-1)2^{m+n+1}})^\top x-(\check{b}_k^{i-(k-1)2^{m+n+1}})^\top y+ \check{c}_k^{i-(k-1)2^{m+n+1}} \mbox{ for } \\ \quad &i=(k-1)2^{m+n+1}+1,
\ldots, k2^{m+n+1} \mbox{ with }  k=1,\ldots,l,
\\g_i(x,y,\mu):=&-c^\top_{k}x-(a_{k}^0+\sum^{s}_{j=1}\check{u}^{j}_{i-(k-1)2^{s}-l2^{m+n+1}}a^j_{k})^\top y +b_{k}^0+\sum^{s}_{j=1}\check{u}^{j}_{i-(k-1)2^{s}-l2^{m+n+1}}b^j_{k}  \mbox{ for}\\ \quad &  i=l2^{m+n+1}+(k-1)2^{s}+1,
\ldots, l2^{m+n+1}+k2^{s} \mbox{ with }  k=1,\ldots,q,\\ g_i(x,y,\mu):=&\mu_{0} \mbox{ for } i=l2^{m+n+1}+q2^s+1,\\ g_i(x,y,\mu):=&\mu_{i-l2^{m+n+1}-q2^s-1} \mbox{ for } i=l2^{m+n+1}+q2^s+2,\ldots,l2^{m+n+1}+q2^s+q+1,\\ g_i(x,y,\mu):=&-\mu_0d_0^\top y-\sum\limits_{k=1}^q(\mu_kb^0_k-\mu_kc_k^\top x+\sum^{s}_{j=1}\mu^j_kb^j_k) \mbox{ for }
i=l2^{m+n+1}+q2^s+q+2,\\
g_i(x,y,\mu):=&\big(\mu_k\gamma_{i-(k-1)s-(l2^{m+n+1}+q2^s+q+2)}\big)^2-\big(\mu_k^{i-(k-1)s-(l2^{m+n+1}+q2^s+q+2)}\big)^2 \mbox{ for } \\ \quad & i=(l2^{m+n+1}+q2^s+q+2)+(k-1)s+1,\ldots, (l2^{m+n+1}+q2^s+q+2)+ks \\ \quad &\mbox{with } k=1,\ldots,q,
\end{align*} and \begin{align*}h_j(x,y,\mu):=& \mu_0d_0^j+\sum\limits_{k=1}^q(\mu_ka^{0j}_k+\sum_{i=1}^s\mu_k^ia^{ij}_k) \mbox{ for } j=1,\ldots,n,\\ h_j(x,y,\mu):=&1- \sum\limits_{k=0}^q(\mu_k)^2-\sum\limits_{k=1}^q\sum\limits_{i=1}^s(\mu_k^i)^2 \mbox{ for }  j=n+1.
\end{align*}
\end{Theorem}
{\bf Proof.} [{\bf Equivalent representation by finite number of constraints}]  Let us first show that  the set \begin{align}\label{1.17}\{(x,y)\in\R^m\times\R^n\mid\tilde a_i^\top x+\tilde b_i^\top y\le \tilde c_i,\; \forall (\tilde a_i, \tilde b_i, \tilde c_i)\in\widetilde U_i,\,i=1,\ldots,l\}\end{align} is equivalent to the following one \begin{align}\label{1.18}\{(x,y)\in\R^m\times\R^n\mid \check{a}_i^{k\top} x+\check{b}_i^{k\top} y- \check{c}_i^{k}\le 0,\; k=1,\ldots, 2^{m+n+1},\,i=1,\ldots,l\},\end{align} where $\{(\check{a}^{k}_i,\check{b}^{k}_i,\check{c}^{k}_i)\in\R^m\times\R^n\times\R\mid k=1,\ldots,2^{m+n+1}\}$ are the  extreme points of the box $\widetilde U_i$ for $i=1,\ldots,l$ as denoted above. Indeed, for each $(x,y)\in\R^m\times\R^n$, by letting $X:=(x,y,-1)$, we obtain that
\begin{align*}&\max{\big\{\tilde a_i^\top x+\tilde b_i^\top y-\tilde c_i\mid (\tilde a_i,\tilde b_i,\tilde c_i)\in
\widetilde U_i, \, i=1,\ldots,l\big\}}\\&=\max{\big\{\tilde a_i^\top x+\tilde b_i^\top y-\tilde c_i\mid (\tilde a_i,\tilde b_i,\tilde c_i)\in {\rm
conv\,}\{(\check{a}^{k}_i,\check{b}^{k}_i,\check{c}^{k}_i), k=1,\ldots,
2^{m+n+1}\},\, i=1,\ldots,l\big\}}\\&=\max{\big\{ W_i^\top
X\mid W_i\in {\rm conv\,}\{(\check{a}^{k}_i,\check{b}^{k}_i,\check{c}^{k}_i), k=1,\ldots,
2^{m+n+1}\}, i=1,\ldots,l\big\}}\\&=\max{\big\{\check{W}_i^{k\top} X\mid k=1,\ldots, 2^{m+n+1},\, i=1,\ldots,l
\big\}}\\&=\max{\big\{\check{a}_i^{k\top} x+\check{b}_i^{k\top}y-\check{c}^{k}_i\mid  k=1,\ldots, 2^{m+n+1},\,i=1,\ldots,l \big\}},
\end{align*}
where $W_i:=(\tilde a_i,\tilde b_i,\tilde c_i)$ and $\check{W}^{k}_i:=(\check{a}^{k}_i,\check{b}^{k}_i,\check{c}^{k}_i), k=1,\ldots,
2^{m+n+1}$ for $i=1,\ldots,l.$ So, we conclude  that  \eqref{1.17} is equivalent to \eqref{1.18}.

\noindent[{\bf Characterizing robust feasible points of~(P)}] Let  $(x,y)\in\R^m\times\R^n$ be a  robust feasible point of problem~(P).  It means that \begin{align}\label{1.19}y\in Y(x),\; \tilde a_i^\top x+\tilde b_i^\top y\le \tilde c_i,\; \forall (\tilde a_i, \tilde b_i, \tilde c_i)\in\widetilde U_i,\; i=1,\ldots,l.\end{align} Due to the equivalence between \eqref{1.17} and \eqref{1.18},  \eqref{1.19} is nothing else but the assertion that
\begin{align}\label{1.16}y\in Y(x),\;  \check{a}_i^{k\top} x+\check{b}_i^{k\top} y- \check{c}_i^{k}\le 0,\; k=1,\ldots, 2^{m+n+1},\; i=1,\ldots,l.\end{align} Invoking Corollary~\ref{Theo2}, \eqref{1.16} can be equivalently expressed as \begin{align*}({\rm I})&\begin{cases}&(x,y)\in\R^m\times\R^n,\\&  \check{a}_i^{k\top} x+\check{b}_i^{k\top} y- \check{c}_i^{k}\le 0,\; k=1,\ldots, 2^{m+n+1},\; i=1,\ldots,l,\\& c^\top_jx +(a_j^0+\sum^{s}_{i=1}\check{u}^{i}_ka^i_j)^\top y -(b_j^0+ \sum^{s}_{i=1}\check{u}^{i}_kb^i_j)\le 0, k=1,\ldots, 2^s, j=1,\ldots,q, \\& \exists \mu_0> 0, \mu_j\ge 0, \mu_j^i\in\R, j=1,\ldots,q, i=1,\ldots,s \mbox{ such that } \sum\limits_{j=0}^q(\mu_j)^2+\sum\limits_{j=1}^q\sum\limits_{i=1}^s(\mu_j^i)^2=1, \\ &\mu_0d_0+\sum\limits_{j=1}^q(\mu_ja^0_j+\sum^{s}_{i=1}\mu^i_ja^i_j)=0, \; -\mu_0d_0^\top y-\sum\limits_{j=1}^q(\mu_jb^0_j-\mu_jc_j^\top x+\sum^{s}_{i=1}\mu^i_jb^i_j)\ge 0 \\&\mbox{and } (\mu_j\gamma_i)^2-(\mu^i_j)^2\ge 0, i=1,\ldots, s,  j=1,\ldots,q,\end{cases} \end{align*}  where $\{\check{u}_k:=(\check{u}^{1}_k,\ldots,\check{u}^{s}_k)\in\R^s\mid k=1,\ldots, 2^s\}$ are the extreme points of the box $V_s.$ Denote $\mu:=(\mu_0,\mu_1,\ldots,\mu_q,\mu_1^1,\ldots,\mu^1_q,\ldots,\mu^s_1,\ldots,\mu^s_q)\in \R^{q(s+1)+1}, g_i(x,y,\mu), i=1,\ldots,L$ and $  h_j(x,y,\mu), j=1,\ldots, n+1$ as stated in the theorem. Then, we conclude by (I) that $(x,y)\in\R^m\times\R^n$ is a robust feasible point of problem~(P) if and only if there exists $\mu\in \R^{q(s+1)+1}$ such that  \begin{align*}({\rm II})\begin{cases} &g_i(x,y,\mu)\ge 0,\; i=1,\ldots, l2^{m+n+1}+q2^s,\\ &g_i(x,y,\mu)>0,\; i=l2^{m+n+1}+q2^s+1,\\ & g_i(x,y,\mu)\ge 0,\; i=l2^{m+n+1}+q2^s+2,\ldots,L,\\
& h_j(x,y,\mu)\ge 0, -h_j(x,y,\mu)\ge 0,\; j=1,\ldots,n+1.\end{cases} \end{align*}

\noindent[{\bf Establishing conditions for applying Putinar's Positivstellensatz}]  We consider a set  of polynomials  $$\mathbf{M}(g_1,\ldots,g_{L},h_1,\ldots,h_{n+1},-h_1,\ldots,-h_{n+1},\kappa-f)$$ as defined in \eqref{Set-polynomial}. It is obvious by  definition  that the polynomial \begin{align*}\hat h:=&\kappa-f+h_{n+1}=1.(\kappa-f)+1.h_{n+1}\\&\in \mathbf{M}(g_1,\ldots,g_{L},h_1,\ldots,h_{n+1},-h_1,\ldots,-h_{n+1},\kappa-f).\end{align*} Since  $(\bar x,\bar y)$ is a  robust feasible point of problem~(P), there exists $\bar\mu\in \R^{q(s+1)+1}$ such that (II) holds at $(\bar x,\bar y,\bar\mu).$ It entails that $h_{n+1}(\bar x,\bar y,\bar\mu)=0$, which in turn implies that   $\hat h(\bar x,\bar y,\bar\mu)=\kappa-f(\bar x,\bar y)\ge 0$, and so,  the set $$H:=\{(x,y,\mu)\in\R^m\times\R^n\times\R^{q(s+1)+1}\mid \hat h(x,y,\mu)\ge 0\}\neq \emptyset$$ by virtue of $(\bar x,\bar y,\bar\mu)\in H$. Take  any $(x,y,\mu)\in H$. We have  $\hat h(x,y,\mu)\ge 0$, which  entails that
\begin{align}\label{1.26}\begin{cases}f(x,y)\le 1+\kappa,\\  \sum\limits_{k=0}^q(\mu_k)^2+\sum\limits_{k=1}^q\sum\limits_{i=1}^s(\mu_k^i)^2 \le 1+\kappa-\inf\limits_{(x,y)\in\R^m\times\R^n}f(x,y).
\end{cases}\end{align} Since $f$ is coercive on $\R^m\times\R^n,$ it follows that  $\inf_{(x,y)\in\R^m\times\R^n}f(x,y)>-\infty$, and hence, \eqref{1.26} guarantees that  $H$ is a compact set. Hence, $$\mathbf{M}(g_1,\ldots,g_{L},h_1,\ldots,h_{n+1},-h_1,\ldots,-h_{n+1},\kappa-f)$$ is archimedean. Put \begin{align*}K:=\{(x,y,\mu)\in \R^m\times\R^n\times\R^{q(s+1)+1}\mid &g_i(x,y,\mu)\ge 0, i=1,\ldots,L, \\&h_j(x,y,\mu)\ge 0, -h_j(x,y,\mu)\ge 0, j=1,\ldots, n+1,\\& \kappa-f(x,y)\ge 0\}.\end{align*} We  easily verify that $(\bar x,\bar y,\bar\mu)\in K,$ and hence, $K\neq\emptyset.$

$[\Longrightarrow]$ Let  $(\bar x,\bar y)$ be  a global  robust solution of problem~(P). For each $\epsilon>0,$ set $\hat f(x,y,\mu):=f(x,y)-f(\bar x,\bar y)+\epsilon$.   We will show that $\hat f>0$ on  $K.$ Indeed, take any $(x,y,\mu)\in K$. It follows that \begin{align*}\begin{cases} &g_i(x,y,\mu)\ge 0,\; i=1,\ldots,L,\\
& h_j(x,y,\mu)\ge 0, -h_j(x,y,\mu)\ge 0,\; j=1,\ldots,n+1,\end{cases} \end{align*} which is equivalent to the following  one \begin{align*}{\rm(III)}&\begin{cases}&(x,y)\in\R^m\times\R^n,\mu:=(\mu_0,\mu_1,\ldots,\mu_q,\mu_1^1,\ldots,\mu^1_q,\ldots,\mu^s_1,\ldots,\mu^s_q)\in \R^{q(s+1)+1},\\&  \check{a}_i^{k\top} x+\check{b}_i^{k\top} y- \check{c}_i^{k}\le 0,\; k=1,\ldots, 2^{m+n+1},\; i=1,\ldots,l,\\& c^\top_jx +(a_j^0+\sum^{s}_{i=1}\check{u}^{i}_ka^i_j)^\top y -(b_j^0+ \sum^{s}_{i=1}\check{u}^{i}_kb^i_j)\le 0, k=1,\ldots, 2^s, j=1,\ldots,q, \\& \mu_0\ge 0, \mu_j\ge 0,  j=1,\ldots,q, \; \sum\limits_{j=0}^q(\mu_j)^2+\sum\limits_{j=1}^q\sum\limits_{i=1}^s(\mu_j^i)^2=1, \\ &\mu_0d_0+\sum\limits_{j=1}^q(\mu_ja^0_j+\sum^{s}_{i=1}\mu^i_ja^i_j)=0, \; -\mu_0d_0^\top y-\sum\limits_{j=1}^q(\mu_jb^0_j-\mu_jc_j^\top x+\sum^{s}_{i=1}\mu^i_jb^i_j)\ge 0 \\&\mbox{and } (\mu_j\gamma_i)^2-(\mu^i_j)^2\ge 0, i=1,\ldots, s,  j=1,\ldots,q.\end{cases} \end{align*} Since the {\rm (LSC)} in \eqref{LLSC} is satisfied, we  will show that $\mu_0\neq 0.$ Assume on the contrary that $\mu_0=0.$ We get by (III) that \begin{align}\label{1.22}&\sum\limits_{j=1}^q(\mu_j)^2+\sum\limits_{j=1}^q\sum\limits_{i=1}^s(\mu_j^i)^2=1, \\ \label{1.23}& \sum\limits_{j=1}^q(\mu_ja^0_j+\sum^{s}_{i=1}\mu^i_ja^i_j)=0, \quad \sum\limits_{j=1}^q \mu_jc_j^\top x -\sum\limits_{j=1}^q\mu_jb^0_j-\sum\limits_{j=1}^q\sum^{s}_{i=1}\mu^i_jb^i_j\ge 0, \\ \label{1.24}& (\mu_j\gamma_i)^2-(\mu^i_j)^2\ge 0, i=1,\ldots, s,  j=1,\ldots,q.\end{align}
By \eqref{1.24}, for each $j\in\{1,\ldots,q\}$, $\mu_j^i=0$ for $i=1,\ldots,s$ whenever $\mu_j=0,$ and hence, we conclude from \eqref{1.22} that there exists $j\in\{1,\ldots,q\}$ such that $\mu_j\neq 0$, i.e., $\sum\limits_{j=1}^q\mu_j\neq 0.$  For each $j\in\{1,\ldots,q\},$ let $\bar u_j:=(\bar u^1_j,\ldots,\bar u^s_j)\in\R^s$ with $$\bar u^i_j:=\begin{cases}0 &\mbox{if } \mu_j=0\\ \frac{\mu^i_j}{\mu_j} &\mbox{if } \mu_j\neq 0,\end{cases}\; i=1,\ldots, s.$$ Then, $|\bar u^i_j|\le \gamma_i$ for $i=1,\ldots,s, j=1,\ldots,q$, and so, $\bar u_j\in V_s=U_j$ for all $j=1,\ldots,q.$ On the one hand, in view of the (LSC) in~\eqref{LLSC}, we find $z\in\R^n$ such that \begin{align*}  c_j^\top x +a_j(\bar u_j)^\top z< b_j(\bar u_j),\; j=1,\ldots,q.\end{align*}

Since  $\sum\limits_{j=1}^q\mu_j\neq 0,$ it follows that  \begin{align*}  \sum\limits_{j=1}^q\mu_j c_j^\top x +\sum\limits_{j=1}^q\mu_j a_j(\bar u_j)^\top z< \sum\limits_{j=1}^q\mu_jb_j(\bar u_j),\end{align*} or equivalently, \begin{align}\label{1.25}  \sum\limits_{j=1}^q\mu_j c_j^\top x +\big(\sum\limits_{j=1}^q\mu_j a_j^0+\sum\limits_{j=1}^q\sum^{s}_{i=1}\mu^i_ja^i_j\big)^\top z- \sum\limits_{j=1}^q\mu_jb_j^0-\sum\limits_{j=1}^q\sum^{s}_{i=1}\mu^i_jb^i_j<0.\end{align} On the other hand, we get by \eqref{1.23} that \begin{align*}  \sum\limits_{j=1}^q\mu_j c_j^\top x +\big(\sum\limits_{j=1}^q\mu_j a_j^0+\sum\limits_{j=1}^q\sum^{s}_{i=1}\mu^i_ja^i_j\big)^\top z- \sum\limits_{j=1}^q\mu_jb_j^0-\sum\limits_{j=1}^q\sum^{s}_{i=1}\mu^i_jb^i_j\ge 0,\end{align*} which contradicts \eqref{1.25}. Therefore, we conclude that $\mu_0> 0.$ It follows that $(x,y,\mu)$ satisfies (II), and so,  $(x,y)$ is a robust feasible point of problem~(P). This fact gives us that $f(x,y)\ge f(\bar x,\bar y)$  inasmuch as $(\bar x,\bar y)$ is  a global  robust solution of problem~(P). Consequently,  it guarantees that $\hat f(x,y,\mu):=f(x,y)-f(\bar x,\bar y)+\epsilon>0.$

\noindent[{\bf Applying Putinar's Positivstellensatz}]
Now, applying Lemma~\ref{Putinar}, we find sums-of-squares polynomials $\sigma_i, i=0,1,\ldots,L, \xi^1_j, \xi^2_j,  j=1,\ldots, n+1, \zeta\in \Sigma^2[x,y,\mu]$ such that
\begin{align}\label{1.27}\hat f=\sigma_0+\sum\limits_{i=1}^{L}\sigma_ig_i+\sum\limits_{j=1}^{n+1}\xi^1_jh_j+\sum\limits_{j=1}^{n+1}\xi^2_j(-h_j)+\zeta(\kappa-f).\end{align} Let $\xi_j\in\R[x,y,\mu],   j=1,\ldots, n+1$ be real polynomials defined by $\xi_j:= \xi^1_j-\xi^2_j,  j=1,\ldots, n+1.$ Then, we deduce from \eqref{1.27} that \begin{align*}\label{}f-\sum\limits_{i=1}^{L}\sigma_ig_i-\sum\limits_{j=1}^{n+1}\xi_jh_j-\zeta(\kappa-f)-f(\bar x,\bar y)+\epsilon=\sigma_0,\end{align*} showing \eqref{1.28} is valid.

$[\Longleftarrow]$ Assume that for any $\epsilon>0$, there exist sums-of-squares polynomials  $\sigma_i, i=0,1,\ldots,L, \zeta\in \Sigma^2[x,y,\mu]$ and real polynomials $\xi_j\in\R[x,y,\mu],   j=1,\ldots, n+1$ such that \eqref{1.28} holds. By rearranging \eqref{1.28}, we obtain that \begin{align}\label{2.9}(1+\zeta)f=\sigma_0+\sum\limits_{i=1}^{L}\sigma_ig_i+\sum\limits_{j=1}^{n+1}\xi_jh_j+(1+\zeta)f(\bar x,\bar y)+\zeta\big(\kappa-f(\bar x,\bar y)\big)-\epsilon.\end{align} Let $(x,y)\in\R^m\times\R^n$ be a  robust feasible point of problem~(P). Then,  there exists $\mu\in \R^{q(s+1)+1}$ such that (II) holds at $(x,y,\mu)$. Thanks to the  nonnegativity of sums-of-squares polynomials,  evaluating  \eqref{2.9} at $(x,y,\mu)$ allows us to arrive at  $$\big(1+\zeta(x,y,\mu)\big)f(x,y)\ge \big(1+\zeta(x,y,\mu)\big)f(\bar x,\bar y)-\epsilon$$ and then, $$f(x,y)\ge f(\bar x,\bar y)- \frac{\epsilon}{1+\zeta(x,y,\mu)}\ge f(\bar x,\bar y)-\epsilon.$$
Letting $\epsilon\to 0$, we obtain that $f(x,y)\ge f(\bar x,\bar y).$ In conclusion, $(\bar x,\bar y)$ is  a global  robust solution of problem~(P), which ends the proof of the theorem. $\hfill\Box$

\begin{Remark}{\rm (i)\ It is worth noticing that if $f$ is a convex polynomial and there exists $(\bar x,\bar y)\in\R^m\times\R^n$ such that the Hessian $\nabla^2f(\bar x,\bar y)$ is positive definite, then it is coercive on $\R^m\times\R^n$ (see e.g., \cite[Lemma~3.1]{jeya-Son-Li-ORL14}) and hence, the above theorem can be obviously applied for this convex setting.\\
(ii)\ Observe  from  the converse implication of the proof of Theorem~\ref{Theo5} that  a  robust feasible point $(\bar x,\bar y)\in\R^m\times\R^n$  is a global robust solution  of problem~(P) if it satisfies  the  following representation  \begin{align}\label{1.28-without-e}f-\sum\limits_{i=1}^{L}\sigma_ig_i-\sum\limits_{j=1}^{n+1}\xi_jh_j-\zeta\big(\kappa-f\big)-f(\bar x,\bar y)=\sigma_0,\end{align} where $\sigma_i, i=0,1,\ldots,L, \zeta\in \Sigma^2[x,y,\mu],$ $\xi_j\in\R[x,y,\mu],   j=1,\ldots, n+1$ and $g_i, i=1,\ldots,L, h_j, j=1,\ldots, n+1$ are defined as in the statement of Theorem~\ref{Theo5}. It is clear that the condition~\eqref{1.28-without-e}  serves as a {\it sufficient criterion}  for global robust optimality, and furthermore, in practice, it is easily  numerically checkable   compared to  the condition~\eqref{1.28},  because it does not involve the parameter $\epsilon.$
}\end{Remark}

The next example illustrates that if the {\rm (LSC)} in \eqref{LLSC} is violated, then the conclusion of  Theorem~\ref{Theo5} may go awry. Below, we consider the case of $l:=1$ and $q:=1$ for the purpose of simplicity.

\begin{Example}\label{Exam1} {\bf (The importance of the Slater condition)} {\rm Consider  the bilevel polynomial  optimization problem with uncertain linear constraints of the form: \begin{align}\label{EP}  \min_{(x,y)\in \R^2}{\big\{f(x,y):=x^2+y^2+2y-2} \mid  y\in Y(x, u)\big\},\tag{EP1} \end{align} where $u\in U:=[-1,1]\subset\R$ and $Y(x, u):={\rm argmin}_{z\in\R}\{x-z\mid (1+u)z\le 0\}.$ The problem~\eqref{EP} can be expressed in terms of problem~(P), where $\widetilde U:=[\underline{a},\overline{a}]\times [\underline{b},\overline{b}]\times[\underline{c},\overline{c}]$ with  $\underline{a}:=\overline{a}:=\underline{b}:=\overline{b}:=\underline{c}:=\overline{c}:=0\in\R$, and $c_0:=1\in\R, d_0:=-1\in\R, c:=0\in\R,$  the affine mappings $a:\R\to \R$ and $b:\R\to\R$ are given by $a(u):=a^0+ua^1$ and $b(u):=b^0+ub^1$ with $ a^0:=a^1:=1\in\R,  b^0:=b^1:=0\in\R$ for $u\in\R$.

In this setting, it is easy to see that $Y(x)=\{0\}$ for any $x\in\R,$ and therefore, we see that $(\bar x,\bar y):=(0,0)$ is a global robust solution of problem~\eqref{EP}. Obviously, $f$ is coercive on $\R^2$.  Let $\mu:=(\mu_0,\mu_1,\mu^1_1)\in\R^3$, and denote the functions $g_i(x,y,\mu)$ and $h_j(x,y,\mu)$ as in the statement of Theorem~\ref{Theo5}. For the sake of clear representation, we first remove the null functions and then relabel them as  \begin{align*}\begin{cases} &g_1(x,y,\mu)=-2y,\;g_2(x,y,\mu)=\mu_0,\; g_3(x,y,\mu)=\mu_1,\\&g_4(x,y,\mu)=\mu_0y,\; g_5(x,y,\mu)=(\mu_1)^2-(\mu_1^1)^2,\\
& h_1(x,y,\mu)=-\mu_0+\mu_1+\mu_1^1,\; h_2(x,y,\mu)=1-(\mu_0)^2-(\mu_1)^2-(\mu_1^1)^2.\end{cases} \end{align*}
Let $\kappa\ge -2= f(\bar x,\bar y)$, and let $0<\epsilon<1$. We claim that the representation of sums-of-squares polynomials given in \eqref{1.28} of Theorem~\ref{Theo5} fails to hold. Indeed, assume on the contrary that there exist sums-of-squares polynomials  $\zeta, \sigma_0, \sigma_i\in \Sigma^2[x,y,\mu], i=1,\ldots,5,$ and real polynomials $\xi_j\in\R[x,y,\mu],   j=1,2$ such that \begin{align}\label{E1.28}f-\sum\limits_{i=1}^{5}\sigma_ig_i-\sum\limits_{j=1}^{2}\xi_jh_j-\zeta\big(\kappa-f\big)-f(\bar x,\bar y)+\epsilon=\sigma_0\end{align}  Setting $\tilde x:=0, \tilde y:=-1$ and $ \tilde\mu:=(0,\frac{1}{\sqrt{2}},\frac{-1}{\sqrt{2}})$, and then substituting $(\tilde x,\tilde y,\tilde \mu)$ into \eqref{E1.28} we obtain that $$-1-2\sigma_1(\tilde x,\tilde y,\tilde \mu)-\frac{1}{\sqrt{2}}\sigma_3(\tilde x,\tilde y,\tilde \mu)-(\kappa+3)\zeta(\tilde x,\tilde y,\tilde \mu)+\epsilon=\sigma_0(\tilde x,\tilde y,\tilde \mu),$$ and hence, it entails that $\sigma_0(\tilde x,\tilde y,\tilde \mu)\le \epsilon-1<0$, which is absurd.

Consequently, the conclusion  of Theorem~\ref{Theo5} fails. The reason is that  the {\rm (LSC)} in \eqref{LLSC} is violated.
} \end{Example}

\section{Global Optimal Values  by Semidefinite Programs}
\setcounter{equation}{0}
We now present semidefinite programming relaxations for the bilevel  polynomial optimization problem with  uncertain linear  constraints~(P) and show how the global optimal value of the bilevel polynomial problem can be found by solving a sequence of corresponding semidefinite programming relaxation problems.

\medskip
For each $k\in\N$, the {\it sums-of-squares optimization} problem associated with the problem~(P) is given by
\begin{align}\label{Dk}  \sup_{(t,\sigma_0,\sigma_i,\xi_j,\zeta)}{\Big\{t} \mid &f-\sum\limits_{i=1}^{L}\sigma_ig_i-\sum\limits_{j=1}^{n+1}\xi_jh_j-\zeta(\kappa-f)-t=\sigma_0,\notag  \\& \notag t\in\R, \zeta,\sigma_0,\sigma_i\in \Sigma^2[x,y,\mu],\xi_j\in\R[x,y,\mu],  \\&\notag  {\rm deg}(\sigma_0)\le k, {\rm deg}(\zeta f)\le k, {\rm deg}(\sigma_ig_i)\le k, {\rm deg}(\xi_jh_j)\le k,\\&  i=1,\ldots,L, j=1,\ldots, n+1\Big\}, \tag{${\rm D_k}$}\end{align} where $\kappa\ge f(\bar x,\bar y)$ with  $(\bar x,\bar y)$ being a robust feasible point of problem~(P),  and  $g_i, i=1,\ldots,L:=l2^{m+n+1}+q(2^s+s+1)+2, h_j, j=1,\ldots, n+1$ are defined as in the statement of Theorem~\ref{Theo5}.

\medskip
It is worth mentioning here that for each fixed $k\in\N$, the problem \eqref{Dk} can be regarded as a sum of squares relaxation problem of the primal one~(P), and more interestingly, it can be reformulated and solved as a semidefinite linear programming problem \cite{Lasserre-book-09}.  Denote the optimal values of~(P) and \eqref{Dk} respectively by ${\rm val}(P)$ and ${\rm val}(D_k)$.
\medskip

In the next theorem we show that, under some additional conditions, the bilevel polynomial optimization problem~(P) has  a global robust solution and the optimal values of the relaxation problem~\eqref{Dk} converge to the optimal value of the   bilevel  polynomial optimization problem~(P) when the degree bound $k$ goes to infinity.

\begin{Theorem}\label{Theo6}{\bf (Computing global optimal value by SDP)} Let $f$ be coercive on $\R^m\times\R^n$.  Assume that the {\rm (LSC)} in \eqref{LLSC} is satisfied. Then, the bilevel polynomial optimization problem~{\rm(P)} has  a global robust solution $(x_0,y_0)$ satisfying   \begin{align}\label{2.5} {\rm val}(D_k)\le {\rm val}(P)=f(x_0,y_0)\; \mbox{ for all }\;  k\in\N. \end{align} Moreover, we have \begin{align}\label{2.3}\lim\limits_{k\to\infty}{\rm val}(D_k)={\rm val}(P).\end{align}
\end{Theorem}
{\bf Proof.} \noindent[{\bf Proving the existence of global robust solutions of~(P)}]  As shown in the proof of Theorem~\ref{Theo5}, we conclude that  $(x,y)\in\R^m\times\R^n$ is a  robust feasible point of problem~(P) if and only if    there exists $\mu\in \R^{q(s+1)+1}$ such that
\begin{align}\label{2.2}\begin{cases} &g_i(x,y,\mu)\ge 0,\; i=1,\ldots, l2^{m+n+1}+q2^s,\\ &g_i(x,y,\mu)>0,\; i=l2^{m+n+1}+q2^s+1,\\ & g_i(x,y,\mu)\ge 0,\; i=l2^{m+n+1}+q2^s+2,\ldots,L,\\
& h_j(x,y,\mu)\ge 0, -h_j(x,y,\mu)\ge 0,\; j=1,\ldots,n+1,\end{cases} \end{align}
 where  $g_i, i=1,\ldots,L, h_j, j=1,\ldots, n+1$ are defined as in the statement of Theorem~\ref{Theo5}.

Let $k\in\N$, and let  $(\bar x,\bar y)\in\R^m\times\R^n$  be a  robust feasible point of problem~(P) as   in the construction of  problem~\eqref{Dk}. Then, one has $\kappa\ge f(\bar x,\bar y),$ and there exists $\bar\mu\in \R^{q(s+1)+1}$ such that \eqref{2.2} holds at $(\bar x,\bar y,\bar\mu)$.  Now, let  $\epsilon>0$ be fixed and  consider the function $\hat f(x,y,\mu):=f(x,y)-f(\bar x,\bar y)+\epsilon$  for $(x,y,\mu)\in \R^m\times\R^n\times\R^{q(s+1)+1}.$ Under the coercivity of $f$, as shown  in the proof of Theorem~\ref{Theo5}, the set $H:=\{(x,y,\mu)\in\R^m\times\R^n\times\R^{q(s+1)+1}\mid \hat h(x,y,\mu)\ge 0\}$ with $\hat h:=(\kappa-f)+h_{n+1}$ is compact, and the set \begin{align*}K:=\{(x,y,\mu)\in \R^m\times\R^n\times\R^{q(s+1)+1}\mid &g_i(x,y,\mu)\ge 0, i=1,\ldots,L, \\&h_j(x,y,\mu)\ge 0, -h_j(x,y,\mu)\ge 0, j=1,\ldots, n+1,\\& \kappa- f(x,y)\ge 0\}\neq\emptyset\end{align*} by virtue of $(\bar x,\bar y,\bar\mu)\in K.$ It can be checked that $K\subset H$, and hence, $K$ is compact as well. In addition, since $\hat f$ is a polynomial and hence continuous, we conclude that there exists $(x_0,y_0,\mu_0)\in K$ such that \begin{align}\label{2.4}\hat f(x_0,y_0,\mu_0)\le \hat f(x,y,\mu)\; \mbox{ for all } (x,y,\mu)\in K.\end{align}
We claim that $(x_0,y_0)$ is a global robust solution of problem~(P). Indeed, by $(x_0,y_0,\mu_0)\in K,$ it follows that \begin{align}\label{2.12}\begin{cases} &g_i(x_0,y_0,\mu_0)\ge 0,\; i=1,\ldots,L,\\
& h_j(x_0,y_0,\mu_0)\ge 0, -h_j(x_0,y_0,\mu_0)\ge 0,\; j=1,\ldots,n+1,\end{cases} \end{align} and that  \begin{align}\label{2.10}\kappa\ge f(x_0,y_0).\end{align} Under the fulfilment of the {\rm (LSC)} in \eqref{LLSC}, similar to the proof of  Theorem~\ref{Theo5}, the condition~\eqref{2.12} guarantees that  \eqref{2.2} holds at $(x_0,y_0,\mu_0)$ and hence, $(x_0,y_0)$ is  a  robust feasible point of problem~(P).

Now, let  $(x,y)\in\R^m\times\R^n$ be a  robust feasible point of problem~(P). Then, there is $\mu\in \R^{q(s+1)+1}$ such that \eqref{2.2} holds. If in addition $ \kappa-f(x,y)\ge 0$, then $(x,y,\mu)\in K$, and so, we get by \eqref{2.4} that $f(x_0,y_0)\le f(x,y).$ Otherwise, $ \kappa-f(x,y)< 0,$ then $f(x,y)> \kappa \ge f(x_0,y_0),$ where the last inequality holds by virtue of \eqref{2.10}. Consequently, our claim holds.

Furthermore, it confirms that \begin{align}\label{2.11}{\rm val}(P)=f(x_0,y_0).\end{align}

\noindent[{\bf Verifying   \eqref{2.5}}] If the  problem~\eqref{Dk} is not feasible, then ${\rm  val}(D_k)=-\infty$, and in this case,  \eqref{2.5} holds trivially. Now, let $(t,\sigma_0,\sigma_i,\xi_j,\zeta), i=1,\ldots,L, j=1,\ldots, n+1 $ be a feasible point of problem~\eqref{Dk}. It means that there exist $ t\in\R, \zeta,\sigma_0,\sigma_i\in \Sigma^2[x,y,\mu],\xi_j\in\R[x,y,\mu],  {\rm deg}(\sigma_0)\le k, {\rm deg}(\zeta f)\le k, {\rm deg}(\sigma_ig_i)\le k, {\rm deg}(\xi_jh_j)\le k, i=1,\ldots,L, j=1,\ldots, n+1$ such that $$f-\sum\limits_{i=1}^{L}\sigma_ig_i-\sum\limits_{j=1}^{n+1}\xi_jh_j-\zeta(\kappa-f)-t=\sigma_0$$ or equivalently, \begin{align}\label{2.1}(1+\zeta)f=\sigma_0+\sum\limits_{i=1}^{L}\sigma_ig_i+\sum\limits_{j=1}^{n+1}\xi_jh_j+t+\zeta f(x_0,y_0)+\zeta\big(\kappa-f(x_0,y_0)\big),\end{align} where $(x_0,y_0)$ is the global robust solution of problem~(P) as shown above.
Recall here that \eqref{2.10} and \eqref{2.2} hold at $(x_0,y_0,\mu_0).$  Due to the  nonnegativity of sums-of-squares polynomials,  estimating  \eqref{2.1} at $(x_0,y_0,\mu_0)$, we obtain that
\begin{align*}\label{}\big(1+\zeta(x_0,y_0,\mu_0)\big)f(x_0,y_0)\ge t+\zeta(x_0,y_0,\mu_0)f(x_0,y_0),\end{align*} or equivalently, \begin{align*}\label{}f(x_0,y_0)\ge t.\end{align*} It confirms that ${\rm val}(D_k)\le f(x_0,y_0)$, which together with \eqref{2.11} proves that \eqref{2.5} is valid.

\noindent[{\bf Verifying   \eqref{2.3}}] Let $\epsilon>0.$ As shown above, $(x_0,y_0)$  is  a global robust solution of problem~(P) satisfying $\kappa\ge f(\bar x,\bar y)\ge f(x_0,y_0).$ Invoking Theorem~\ref{Theo5}, we find sums-of-squares polynomials $\sigma_i, i=0,1,\ldots,L, \xi^1_j, \xi^2_j,  j=1,\ldots, n+1, \zeta\in \Sigma^2[x,y,\mu]$ such that \begin{align*}\label{}f-\sum\limits_{i=1}^{L}\sigma_ig_i-\sum\limits_{j=1}^{n+1}\xi_jh_j-\zeta(\kappa-f)-f(x_0,y_0)+\epsilon=\sigma_0,\end{align*} or equivalently, \begin{align*}\label{}f-\sum\limits_{i=1}^{L}\sigma_ig_i-\sum\limits_{j=1}^{n+1}\xi_jh_j-\zeta(\kappa-f)-t_0=\sigma_0,\end{align*}
with $t_0:=f(x_0,y_0)-\epsilon\in\R.$ So, there exists $k_\epsilon\in\N$ such that
${\rm deg}(\sigma_0)\le k_\epsilon, {\rm deg}(\zeta f)\le k_\epsilon, {\rm deg}(\sigma_ig_i)\le k_\epsilon, {\rm deg}(\xi_jh_j)\le k_\epsilon,  i=1,\ldots,L, j=1,\ldots, n+1,$ and that $${\rm val}(D_{k_\epsilon})\ge f(x_0,y_0)-\epsilon.$$ Letting $\epsilon\to 0$, we see that $\liminf\limits_{k\to\infty}{\rm val}(D_k)\ge f(x_0,y_0).$ This together with \eqref{2.5} establishes \eqref{2.3}, which ends the proof of the theorem. $\hfill\Box$

\medskip
Finally, we provide some examples which show  how our  relaxation scheme can be applied to find the global optimal value of the bilevel polynomial optimization problem with  uncertain linear  constraints~(P).

\begin{Example}\label{Exam2} {\bf(Uncertain bilevel convex polynomial problem)}  {\rm Consider  the bilevel polynomial  optimization problem with uncertain linear constraints of the form: \begin{align}\label{EP2}  \min_{(x,y)\in \R^2}{\big\{f(x,y):=x^4+y^2+y+1} \mid  y\in Y(x, u), x\le 0, y\le 0\big\},\tag{EP2} \end{align} where $u\in U:=[-1,1]\subset\R$ and $Y(x, u):={\rm argmin}_{z\in\R}\{2x-z\mid (1+\frac{1}{2}u)z\le 0\}.$ The problem~\eqref{EP2} can be expressed in terms of problem~(P), where $\widetilde U:=[\underline{a},\overline{a}]\times [\underline{b},\overline{b}]\times[\underline{c},\overline{c}]$ with  $\underline{a}:=\underline{b}:=\underline{c}:=0\in\R,$ $\overline{a}:=\overline{b}:=\overline{c}:=1\in\R$, and $c_0:=2\in\R, d_0:=-1\in\R, c:=0\in\R,$  the affine mappings $a:\R\to \R$ and $b:\R\to\R$ are given by $a(u):=a^0+ua^1$ and $b(u):=b^0+ub^1$ with $ a^0:=1\in\R, a^1:=\frac{1}{2}\in\R,  b^0:=b^1:=0\in\R$ for $u\in\R$. A  direct calculation shows that $(x_0,y_0):=(0,0)$ is a robust solution of problem~\eqref{EP2} with the global optimal value $1.$

Now, we use the relaxation scheme formulated in Theorem~\ref{Theo6} to verify this global optimal value.  In this setting, it is easy to see that $f$ is coercive on $\R^2$ and the  {\rm (LSC)} in \eqref{LLSC} is fulfilled.  Let $\mu:=(\mu_0,\mu_1,\mu^1_1)\in\R^3$, and denote the functions $g_i(x,y,\mu)$ and $h_j(x,y,\mu)$ as in the statement of Theorem~\ref{Theo5}. For the sake of clear representation, we first remove the null functions and then relabel them as  \begin{align*}\begin{cases} &g_1(x,y,\mu)=1,\;g_2(x,y,\mu)=-x,\;g_3(x,y,\mu)=-x+1,\; g_4(x,y,\mu)=-y,\\&g_5(x,y,\mu)=-y+1,\;g_6(x,y,\mu)=-x-y,\;g_7(x,y,\mu)=-x-y+1,\;g_8(x,y,\mu)=\mu_0,\\&g_9(x,y,\mu)=\mu_1,\; g_{10}(x,y,\mu)=\mu_0y,\; g_{11}(x,y,\mu)=(\mu_1)^2-(\mu_1^1)^2,\\
& h_1(x,y,\mu)=-\mu_0+\mu_1+\frac{1}{2}\mu_1^1,\; h_2(x,y,\mu)=1-(\mu_0)^2-(\mu_1)^2-(\mu_1^1)^2.\end{cases} \end{align*}
Let $(\bar x,\bar y):=(-1,0)$ be a feasible point, and take $\kappa:=2\ge 2= f(\bar x,\bar y)$. In this setting, the   problem~\eqref{Dk} becomes  \begin{align*}\label{}  \sup_{(t,\sigma_0,\sigma_i,\xi_j,\zeta)}{\Big\{t} \mid &f-\sum\limits_{i=1}^{11}\sigma_ig_i-\sum\limits_{j=1}^{2}\xi_jh_j-\zeta(2-f)-t=\sigma_0,\notag  \\& \notag t\in\R, \zeta,\sigma_0,\sigma_i\in \Sigma^2[x,y,\mu],\xi_j\in\R[x,y,\mu],  {\rm deg}(\sigma_0)\le k,\\& {\rm deg}(\zeta f)\le k, {\rm deg}(\sigma_ig_i)\le k, {\rm deg}(\xi_jh_j)\le k, i=1,\ldots,11, j=1,2\Big\}.\end{align*} Using the Matlab toolbox YALMIP \cite{Lof-09,Lof-04}, we converted the above optimization problem into an equivalent semi-definite program and solved it  with  $k:=6$. The solver returned the true global optimal value $1.000$.
} \end{Example}

\begin{Example}\label{Exam3} {\bf(Uncertain bilevel non-convex polynomial problem)} {\rm
Consider  the bilevel polynomial  optimization problem with uncertain linear constraints of the form: 
\begin{align}\label{EP3}  \min_{(x,y)\in \R^2}{\big\{f(x,y):=x^4-4xy+y^4-2} \mid  y\in Y(x, u)\big\},\tag{EP3} \end{align} where $u\in U:=[-\frac{1}{2},\frac{1}{2}]\subset\R$ and $Y(x, u):={\rm argmin}_{z\in\R}\{x-z\mid (1+u)z\le 0\}.$ 
The problem~\eqref{EP3} can be expressed in terms of problem~(P), where $\widetilde U:=[\underline{a},\overline{a}]\times [\underline{b},\overline{b}]\times[\underline{c},\overline{c}]$ with  $\underline{a}:=\overline{a}:=\underline{b}:=\overline{b}:=\underline{c}:=\overline{c}:=0\in\R$, and $c_0:=1\in\R, d_0:=-1\in\R, c:=0\in\R,$  the affine mappings $a:\R\to \R$ and $b:\R\to\R$ are given by $a(u):=a^0+ua^1$ and $b(u):=b^0+ub^1$ with $ a^0:=a^1:=1\in\R,  b^0:=b^1:=0\in\R$ for $u\in\R$. Note that the objective function $f$ is a non-convex polynomial. A  direct calculation shows that $(x_0,y_0):=(0,0)$ is a robust solution of problem~\eqref{EP3} with the global optimal value $-2.$

Now, we employ the relaxation scheme of Theorem~\ref{Theo6} to verify this global optimal value. It is easy to check that  $f$ is coercive on $\R^2$ and  {\rm (LSC)} in \eqref{LLSC} is fulfilled.  Let $\mu:=(\mu_0,\mu_1,\mu^1_1)\in\R^3$, and denote the functions $g_i(x,y,\mu)$ and $h_j(x,y,\mu)$ as in the statement of Theorem~\ref{Theo5}. Removing the null functions and then relabeling them gives us
\begin{align*}\begin{cases} &g_1(x,y,\mu)=-\frac{1}{2}y,\;g_2(x,y,\mu)=-\frac{3}{2}y,\;g_3(x,y,\mu)=\mu_0,\\& g_4(x,y,\mu)=\mu_1,\; g_5(x,y,\mu)=\mu_0y,\; g_6(x,y,\mu)=\frac{1}{4}(\mu_1)^2-(\mu_1^1)^2,\\
& h_1(x,y,\mu)=-\mu_0+\mu_1+\mu_1^1,\; h_2(x,y,\mu)=1-(\mu_0)^2-(\mu_1)^2-(\mu_1^1)^2.\end{cases} \end{align*}
Let $(\bar x,\bar y):=(1,0)$ be a feasible point, and take $\kappa:=-1\ge -1= f(\bar x,\bar y)$. Then,   problem~\eqref{Dk} becomes  \begin{align*}\label{}  \sup_{(t,\sigma_0,\sigma_i,\xi_j,\zeta)}{\Big\{t} \mid &f-\sum\limits_{i=1}^{6}\sigma_ig_i-\sum\limits_{j=1}^{2}\xi_jh_j-\zeta(-1-f)-t=\sigma_0,\notag  \\& \notag t\in\R, \zeta,\sigma_0,\sigma_i\in \Sigma^2[x,y,\mu],\xi_j\in\R[x,y,\mu],  {\rm deg}(\sigma_0)\le k,\\& {\rm deg}(\zeta f)\le k, {\rm deg}(\sigma_ig_i)\le k, {\rm deg}(\xi_jh_j)\le k, i=1,\ldots,6, j=1,2\Big\}.\end{align*} Using the Matlab toolbox YALMIP \cite{Lof-09,Lof-04}, we converted the above optimization problem into an equivalent semi-definite program and solved it  with $k:=4$. The solver returned the true global optimal value $-2.000$. }\end{Example}

\section{Appendix: Bilevel Problems \& Ball Data Uncertainty}
\setcounter{equation}{0}

In this Section, we show how a numerically checkabale characterization of robust feasibility can be derived for a bilevel polynomial optimization problem with uncertain linear constraints in the case of ball data uncertainty.

Let $f:\R^m\times\R^n\rightarrow\R$ be a real polynomial. We consider  a    bilevel polynomial optimization problem with {\it ball uncertainties} as
\begin{align}\label{BP}  \min_{(x,y)\in \R^m\times\R^n}{\big\{f(x,y)} \mid  y\in Y(x, u_1,\ldots, u_q),\;\tilde a_i(\tilde u_i)^\top x\le \tilde b_i(\tilde u_i),\; i=1,\ldots,l\big\}, \tag{BP}\end{align} where $u_j\in U_j, j=1,\ldots, q$ and $\tilde u_i\in \widetilde U_i, i=1,\ldots,l$ are {\it uncertain}
and $Y(x, u_1,\ldots, u_q):={\rm argmin}_{z\in\R^n}\{c_0^\top x+d_0^\top z\mid c_j^\top x +a_j(u_j)^\top z\le b_j(u_j),\;  j=1,\ldots,q\}$ denotes the optimal solution set of the uncertain lower-level optimization problem \begin{align}\label{LBP}\min_{z\in\R^n}{\{c_0^\top x+d_0^\top z\mid c_j^\top x +a_j(u_j)^\top z\le b_j(u_j),\; j=1,\ldots,q\}}.\tag{LBP}\end{align} In the above data, the  {\it uncertainty} sets $U_j:=\B_s, j=1,\ldots,q$ and $\widetilde U_i:=\B_{\tilde s}, i=1,\ldots,l$, where $\B_s$ and $\B_{\tilde s}$ are the closed unit balls in $\R^s$ and $\R^{\tilde s},$ respectively, as well as  $c_0\in\R^m, d_0\in\R^n, c_j \in \R^m,    j=1,\ldots,q$ fixed.
\medskip

The affine mappings $\tilde a_i:\R^{\tilde s}\to \R^n, \tilde b_i:\R^{\tilde s}\to\R, i=1,\ldots,l, a_j:\R^s\to \R^n, b_j:\R^s\to\R, j=1,\ldots,q$ are  given respectively by $\tilde a_i(\tilde u_i):=\tilde a^0_i+\sum^{\tilde s}_{j=1}u^j_ia^j_i, \tilde b_i(\tilde u_i)=\tilde b^0_i+\sum^{\tilde s}_{j=1}u^j_ib^j_i $ for $ \tilde u_i:=(u^1_i,\ldots, u^{\tilde s}_i)\in \R^{\tilde s}$ with $\tilde a_i^j\in \R^n, \tilde b_i^j\in\R, j=0,1,\ldots,\tilde s, i=1,\ldots,l$ fixed, and $a_j(u_j):=a^0_j+\sum^{s}_{i=1}u^i_ja^i_j, b_j(u_j)=b^0_j+\sum^{s}_{i=1}u^i_jb^i_j$  for  $u_j:=(u^1_j,\ldots, u^s_j)\in \R^s$  with $a_j^i\in \R^n, b_j^i\in\R, i=0,1,\ldots,s, j=1,\ldots,q$ fixed.

\medskip
Observe that for each $ j\in\{1,\ldots, q\},$ let \begin{align}\label{matric-ball}A^0_j:= \left(
\begin{array}{cc}
I_s & 0 \\
0& 1 \\
\end{array}
\right),\quad A^i_j:= \left(
\begin{array}{cc}
0 & e_i \\
e^\top_i& 0 \\
\end{array}
\right),\; i=1,\ldots, s,\end{align} where $e_i\in \R^s$ is a vector whose $i$th component is one and all others are zero. Then, we have \begin{align}\label{U-ball}&\big\{u_j:=(u^1_j,\ldots, u^s_j)\in\R^s\mid A^0_j+\sum_{i=1}^su^i_jA^i_j\succeq 0\big\}\notag\\&=\big\{u_j:=(u^1_j,\ldots, u^s_j)\in\R^s\mid \left(
\begin{array}{cc}
I_s & u_j \\
u_j^\top& 1 \\
\end{array}
\right)\succeq 0\big\}\\&=\big\{u_j:=(u^1_j,\ldots, u^s_j)\in\R^s\mid 1-u_j^\top I_su_j\ge 0\big\}\notag\\&=\B_s,\quad j=1,\ldots, q,\notag\end{align} where the third equality in \eqref{U-ball} is valid due to the Schur complement (cf.~\cite[Lemma~4.2.1]{Ben-tal01}). It means that the closed unit ball of $\R^s$ can be expressed in terms of  spectrahedra in \eqref{U-set}.

\medskip  We should note here that the notions of robust feasible/or solutions of the upper/lower-level optimization problems~\eqref{BP} and \eqref{LBP} are defined similarly in the previous sections. The following theorem provides a  characterization for robust solutions of the lower-level optimization problem~\eqref{LBP}.

\begin{Theorem}\label{Theo-App1}{\bf (Characterization for robust solutions of \eqref{LBP})} Let $x\in\R^m$, and let  the cone $C(x):={\rm cone}\big\{\big(a_j(u_j),b_j(u_j)-c_j^\top x\big)\mid u_j\in \B_s, j=1,\ldots,q\big\}$  be closed. Then,  $y\in Y(x)$ if and only if \begin{align*}\begin{cases}&y\in\R^n, \exists\lambda_j\ge 0, j=1,\ldots,q \mbox{ such that }\\&\left(\begin{array}{cc}
\lambda_j I_s &  \frac{1}{2}\left(b^1_j-a^{1\top}_jy,\ldots,b^s_j-a_j^{s\top}y\right)^\top  \\
\frac{1}{2}\left(b^1_j-a^{1\top}_jy,\ldots,b^s_j-a_j^{s\top}y\right) & -\lambda_j-a^{0\top}_jy-c^\top_jx+b_j^0 \\
\end{array}
\right)\succeq 0, j=1,\ldots,q, \\ &\exists  \lambda_j^0\ge 0, \lambda_j^i\in\R, j=1,\ldots,q, i=1,\ldots,s \mbox{ such that }\\&d_0+\sum\limits_{j=1}^q(\lambda_j^0a^0_j+\sum^{s}_{i=1}\lambda^i_ja^i_j)=0, \; -d_0^\top y-\sum\limits_{j=1}^q(\lambda_j^0b^0_j-\lambda_j^0c_j^\top x+\sum^{s}_{i=1}\lambda^i_jb^i_j)\ge 0 \\&\mbox{and } (\lambda_j^0)^2-\sum_{i=1}^s(\lambda^i_j)^2\ge 0,  j=1,\ldots,q,\end{cases}\end{align*} where $I_s$ denotes the identity  $(s\times s)$ matrix.
\end{Theorem}
{\bf Proof.} First, we assert that $y\in \R^n$ is a robust feasible point of  the problem~\eqref{LBP} if and only if there are $\lambda_j\ge 0, j=1,\ldots,q$ such that  \begin{align}\label{1.10}\left(\begin{array}{cc}
\lambda_j I_s &  \frac{1}{2}\left(b^1_j-a^{1\top}_jy,\ldots,b^s_j-a_j^{s\top}y\right)^\top  \\
\frac{1}{2}\left(b^1_j-a^{1\top}_jy,\ldots,b^s_j-a_j^{s\top}y\right) & -\lambda_j-a^{0\top}_jy-c^\top_jx+b_j^0 \\
\end{array}
\right)\succeq 0,\; j=1,\ldots,q. \end{align} Indeed, let $y\in \R^n$ be a robust feasible point of the problem~\eqref{LBP}, i.e., we have  \begin{align}\label{1.9}c_j^\top x +a_j(u_j)^\top y\le b_j(u_j),\;\forall u_j\in \B_s,\; j=1,\ldots,q.\end{align}  The relations in \eqref{1.9} means that   for each $j\in\{1,\ldots,q\},$ the following implication holds: \begin{align*}&\forall u_j:=(u^1_j,\ldots,u^s_j)\in \R^s,\;  f_j(u_j):=u_j^\top I_s u_j-1\le 0 \\ \Longrightarrow &g_j(u_j):= \big[(a^{1\top}_jy,\ldots, a_j^{s\top}y)-(b^1_j,\ldots,b^s_j)\big]^\top u_j+a^{0\top}_jy+c^\top_jx-b_j^0 \le 0.\end{align*}  Moreover, since $f_j(0)=-1<0,$ by using the inhomogeneous $S$-lemma (cf.~\cite[Proposition~4.10.1]{Ben-tal01}), we find $\lambda_j\ge 0$ such that \begin{align*}\label{}\left(\begin{array}{cc}
\lambda_j I_s &  \frac{1}{2}\left(b^1_j-a^{1\top}_jy,\ldots,b^s_j-a_j^{s\top}y\right)^\top  \\
\frac{1}{2}\left(b^1_j-a^{1\top}_jy,\ldots,b^s_j-a_j^{s\top}y\right) & -\lambda_j-a^{0\top}_jy-c^\top_jx+b_j^0 \\
\end{array}\right)\succeq 0, \end{align*} and so, our assertion holds.

Keeping in mind the above fact, due to the closed cone $C(x)$, we  apply Theorem~\ref{Theo2-LN} to conclude that  $y\in Y(x)$ if and only if \begin{align*}\begin{cases}&y\in\R^n, \exists\lambda_j\ge 0, j=1,\ldots,q \mbox{ such that }\\&\left(\begin{array}{cc}
\lambda_j I_s &  \frac{1}{2}\left(b^1_j-a^{1\top}_jy,\ldots,b^s_j-a_j^{s\top}y\right)^\top  \\
\frac{1}{2}\left(b^1_j-a^{1\top}_jy,\ldots,b^s_j-a_j^{s\top}y\right) & -\lambda_j-a^{0\top}_jy-c^\top_jx+b_j^0 \\
\end{array}
\right)\succeq 0, j=1,\ldots,q, \\ &\exists  \lambda_j^0\ge 0, \lambda_j^i\in\R, j=1,\ldots,q, i=1,\ldots,s \mbox{ such that }\\&d_0+\sum\limits_{j=1}^q(\lambda_j^0a^0_j+\sum^{s}_{i=1}\lambda^i_ja^i_j)=0, \; -d_0^\top y-\sum\limits_{j=1}^q(\lambda_j^0b^0_j-\lambda_j^0c_j^\top x+\sum^{s}_{i=1}\lambda^i_jb^i_j)\ge 0 \\&\mbox{and } \lambda_j^0A^0_j+\sum_{i=1}^s\lambda^i_jA^i_j\succeq 0,  j=1,\ldots,q.\end{cases}\end{align*} To complete the proof of the theorem, it remains to prove that the matrix inequalities  $ \lambda_j^0A^0_j+\sum_{i=1}^s\lambda^i_jA^i_j\succeq 0, j=1,\ldots,q$ reduce to \begin{align}\label{1.11}(\lambda_j^0)^2-\sum_{i=1}^s(\lambda^i_j)^2\ge 0, j=1,\ldots,q\end{align}  under our setting. Indeed, let $j\in\{1,\ldots,q\}$ and consider the matrices $A^i_j, i=0,1,\ldots, s$ as given in \eqref{matric-ball}. Then, it can be checked that \begin{align} \label{1.12} \lambda_j^0A^0_j+\sum_{i=1}^s\lambda^i_jA^i_j\succeq 0 \Leftrightarrow \left(
\begin{array}{cc}
\lambda_j^0I_s & \bar u_j \\
\bar u_j^\top& \lambda_j^0 \\
\end{array}
\right)\succeq 0\end{align} where $\bar u_j:=(\lambda^1_j,\ldots, \lambda^s_j)\in\R^s$.

If $\lambda_j^0=0,$ then $\lambda_j^i=0$ for all $i=1,\ldots,s$ as shown in the proof of Theorem~\ref{Theo1} due to the boundedness of the $U_j:=\B_s.$ Then, \eqref{1.11} and \eqref{1.12} are  trivially equivalent.  If $\lambda_j^0\neq 0$, then,  by using the Schur complement (cf.~\cite[Lemma~4.2.1]{Ben-tal01}), \eqref{1.12} amounts to the following one $$\lambda_j^0-\bar u_j^\top (\lambda_j^0I_s)^{-1}\bar u_j\ge 0.$$ Now, it is clear  that \begin{align*}\label{}\lambda_j^0-\bar u_j^\top (\lambda_j^0I_s)^{-1}\bar u_j\ge 0\Leftrightarrow\lambda_j^0-\frac{1}{\lambda_j^0}\left(\sum_{i=1}^s(\lambda^i_j)^2\right)\ge 0\Leftrightarrow (\lambda_j^0)^2-\sum_{i=1}^s(\lambda^i_j)^2\ge 0,\end{align*} and so, we arrive at the desired conclusion. $\hfill\Box$

\medskip
The next theorem presents a  characterizing for robust feasible points of the bilevel polynomial optimization problem with  ball uncertainties~\eqref{BP}.

\begin{Theorem}\label{Theo-App2} {\bf (Characterization for robust feasible points of \eqref{BP})}  Let  the cone $C(x):={\rm cone}\big\{\big(a_j(u_j),b_j(u_j)-c_j^\top x\big)\mid u_j\in \B_s, j=1,\ldots,q\big\}$  be closed for each $x\in\R^m.$ Then,  $(x,y)\in\R^m\times\R^n$ is  a  robust feasible point of problem~\eqref{BP} if and only if there exists $\lambda:=(\tilde\lambda_1,\ldots,\tilde\lambda_l,\lambda_1,\ldots,\lambda_q,\lambda^0_1,\ldots,\lambda^0_q,\lambda_1^1,\ldots,\lambda^1_q,\ldots,\lambda^s_1,\ldots,
\lambda^s_q)\in\R^{q(s+2)+l}$ such that \begin{align*}({\rm I})\begin{cases}&\tilde\lambda_i\ge 0, i=1,\ldots,l,\\& \left(\begin{array}{cc}
\tilde\lambda_i I_{\tilde s} &  \frac{1}{2}\left(\tilde b^1_i-\tilde a^{1\top}_ix,\ldots,\tilde b^{\tilde s}_i-\tilde a_i^{{\tilde s}\top}x\right)^\top  \\
\frac{1}{2}\left(\tilde b^1_i-\tilde a^{1\top}_ix,\ldots,\tilde b^{\tilde s}_i-\tilde a_i^{{\tilde s}\top}x\right) & -\tilde\lambda_i-\tilde a^{0\top}_ix+\tilde b_i^0 \\
\end{array}
\right)\succeq 0, i=1,\ldots,l, \\
&\lambda_j\ge 0, j=1,\ldots,q,\\&\left(\begin{array}{cc}
\lambda_j I_s &  \frac{1}{2}\left(b^1_j-a^{1\top}_jy,\ldots,b^s_j-a_j^{s\top}y\right)^\top  \\
\frac{1}{2}\left(b^1_j-a^{1\top}_jy,\ldots,b^s_j-a_j^{s\top}y\right) & -\lambda_j-a^{0\top}_jy-c^\top_jx+b_j^0 \\
\end{array}
\right)\succeq 0, j=1,\ldots,q, \\ &\lambda_j^0\ge 0,  j=1,\ldots,q, \\&d_0+\sum\limits_{j=1}^q(\lambda_j^0a^0_j+\sum^{s}_{i=1}\lambda^i_ja^i_j)=0, \; -d_0^\top y-\sum\limits_{j=1}^q(\lambda_j^0b^0_j-\lambda_j^0c_j^\top x+\sum^{s}_{i=1}\lambda^i_jb^i_j)\ge 0, \\& (\lambda_j^0)^2-\sum_{i=1}^s(\lambda^i_j)^2\ge 0,  j=1,\ldots,q,\end{cases}\end{align*} where $I_k$ denotes the identity  $(k\times k)$ matrix for $k\in\N.$
\end{Theorem}
{\bf Proof.} Let  $(x,y)\in\R^m\times\R^n$ be a  robust feasible point of problem~\eqref{BP}.  It means that \begin{align}\label{1.19-app} y\in Y(x),\;\tilde a_i(\tilde u_i)^\top x\le \tilde b_i(\tilde u_i),\; \forall \tilde u_i\in\B_{\tilde s},\; i=1,\ldots,l.\end{align} Similar to the proof of Theorem~\ref{Theo-App1}, the inequalities $$\tilde a_i(\tilde u_i)^\top x\le \tilde b_i(\tilde u_i),\; \forall \tilde u_i\in\B_{\tilde s},\; i=1,\ldots,l$$ are equivalent  to the assertion that there exist  $\tilde \lambda_i\ge 0, i=1,\ldots,l$ such that  \begin{align*}\label{}\left(\begin{array}{cc}
\tilde\lambda_i I_{\tilde s} &  \frac{1}{2}\left(\tilde b^1_i-\tilde a^{1\top}_ix,\ldots,\tilde b^{\tilde s}_i-\tilde a_i^{{\tilde s}\top}x\right)^\top  \\
\frac{1}{2}\left(\tilde b^1_i-\tilde a^{1\top}_ix,\ldots,\tilde b^{\tilde s}_i-\tilde a_i^{{\tilde s}\top}x\right) & -\tilde\lambda_i-\tilde a^{0\top}_ix+\tilde b_i^0 \\
\end{array}
\right)\succeq 0,\; i=1,\ldots,l. \end{align*}

Now, invoking Theorem~\ref{Theo-App1}, \eqref{1.19-app} can be equivalently expressed as \begin{align*}(\rm {II})&\begin{cases}&(x,y)\in\R^m\times\R^n,\\&\exists\tilde\lambda_i\ge 0, i=1,\ldots,l \mbox{ such that}\\& \left(\begin{array}{cc}
\tilde\lambda_i I_{\tilde s} &  \frac{1}{2}\left(\tilde b^1_i-\tilde a^{1\top}_ix,\ldots,\tilde b^{\tilde s}_i-\tilde a_i^{{\tilde s}\top}x\right)^\top  \\
\frac{1}{2}\left(\tilde b^1_i-\tilde a^{1\top}_ix,\ldots,\tilde b^{\tilde s}_i-\tilde a_i^{{\tilde s}\top}x\right) & -\tilde\lambda_i-\tilde a^{0\top}_ix+\tilde b_i^0 \\
\end{array}
\right)\succeq 0,\; i=1,\ldots,l,\\&\exists\lambda_j\ge 0, j=1,\ldots,q \mbox{ such that}\\&\left(\begin{array}{cc}
\lambda_j I_s &  \frac{1}{2}\left(b^1_j-a^{1\top}_jy,\ldots,b^s_j-a_j^{s\top}y\right)^\top  \\
\frac{1}{2}\left(b^1_j-a^{1\top}_jy,\ldots,b^s_j-a_j^{s\top}y\right) & -\lambda_j-a^{0\top}_jy-c^\top_jx+b_j^0 \\
\end{array}
\right)\succeq 0, j=1,\ldots,q, \\ &\exists  \lambda_j^0\ge 0, \lambda_j^i\in\R, j=1,\ldots,q, i=1,\ldots,s \mbox{ such that }\\&d_0+\sum\limits_{j=1}^q(\lambda_j^0a^0_j+\sum^{s}_{i=1}\lambda^i_ja^i_j)=0, \; -d_0^\top y-\sum\limits_{j=1}^q(\lambda_j^0b^0_j-\lambda_j^0c_j^\top x+\sum^{s}_{i=1}\lambda^i_jb^i_j)\ge 0 \\&\mbox{and } (\lambda_j^0)^2-\sum_{i=1}^s(\lambda^i_j)^2\ge 0,  j=1,\ldots,q.\end{cases}\end{align*}
By denoting $\lambda:=(\tilde\lambda_1,\ldots,\tilde\lambda_l,\lambda_1,\ldots,\lambda_q,\lambda^0_1,\ldots,\lambda^0_q,\lambda_1^1,\ldots,\lambda^1_q,\ldots,\lambda^s_1,\ldots,
\lambda^s_q)\in\R^{q(s+2)+l}$,  the conclusion  follows from (II). $\hfill\Box$

\medskip
\begin{Remark}{\rm Under assumptions of Theorem~\ref{Theo-App2}, let $$X:=\{(x,y,\lambda)\in \R^m\times\R^n\times\R^{q(s+2)+l}\mid (x,y,\lambda)\; \mbox{ satisfying } \; {\rm (I)}\}$$ Then,  the bilevel polynomial optimization problem with  ball uncertainties~\eqref{BP} can be converted into a single-level  polynomial program of the form: \begin{align}\label{M}  \min_{(x,y,\mu)\in \R^m\times\R^n\times\R^{q(s+2)+l}}{\big\{f(x,y)} \mid  (x,y,\lambda)\in X\big\}. \tag{M}\end{align} The  problem~\eqref{M} is a   polynomial program  with linear matrix inequality constraints, which has been studied intensively in the literature; see e.g., \cite{Ben-tal01,Lasserre-book-09}.
}\end{Remark}

\noindent\textbf{Acknowledgments.} The authors are grateful to Dr Guoyin Li, University of New South Wales, for his thoughtful discussions and his help in the computer implementation of our methods.

\end{document}